\documentclass[12pt]{article}

\usepackage[margin=1in]{geometry}
\usepackage{hyperref,subfigure}
\usepackage{tikz}
\usetikzlibrary{shapes,arrows}
\usepackage{paralist,amsmath,amsthm,amssymb}

\newcommand{\sr}{\stackrel}
\newcommand{\ov}{\overline}
\newcommand{\ga}{\gamma}
\newcommand{\al}{\alpha}
\newcommand{\be}{\beta}

\newcommand{\la}{\lambda}
\newcommand{\ty}{type}

\numberwithin{equation}{section}
\newtheorem{thm}[equation]{Theorem}

\newtheorem{defn}[equation]{Definition}
\newtheorem{lem}[equation]{Lemma}
\newtheorem{prop}[equation]{Proposition}
\newtheorem{cor}[equation]{Corollary}
\newtheorem{exm}[equation]{Example}

\newtheorem{rem}[equation]{Remark}

\tikzstyle{wv}=[circle,draw=black!90,fill=white!20,thick,inner sep=2pt,minimum width=5pt] \tikzstyle{bv}=[circle,draw=black!90,fill=black!100,thick,inner sep=2pt,minimum width=5pt] 
\tikzstyle{gv}=[circle,draw=black!70,fill=gray!100,thick,inner sep=2pt,minimum width=5pt]

\title{Direct bijective computation of the generating series for 2 and 3-connection coefficients of the symmetric group}

\author{Alejandro H. Morales\\[-0.8ex] 
\small LaCIM\\[-0.8ex]
\small Universit\'e du Qu\'ebec \`a Montr\'eal\\[-0.8ex] 
\small Montr\'eal PQ H3C 3P8, Canada\\
\small\tt ahmorales@lacim.ca\\
 \and Ekaterina A. Vassilieva\\[-0.8ex] 
\small LIX\\[-0.8ex] 
\small Ecole Polytechnique\\[-0.8ex] 
\small 91128 Palaiseau Cedex, France\\
\small\tt katya@lix.polytechnique.fr}

\begin{document}
\maketitle
%\nodate

\begin{abstract}
We evaluate combinatorially certain connection coefficients of the symmetric group that count the number of factorizations of a long cycle as a product of three permutations. Such factorizations admit an important topological interpretation in terms of unicellular constellations on orientable surfaces.
Algebraic computation of these coefficients was first done by Jackson using irreducible characters of the symmetric group. However, bijective computations of these coefficients are so far limited to very special cases.
% computing these coefficients leads to very intricate formulae with unclear combinatorial interpretation. 
Thanks to a new bijection that refines the work of Schaeffer and Vassilieva in \cite{SV} and Vassilieva in \cite{V}, we give an explicit closed form evaluation of the generating series for these coefficients. The main ingredient in the bijection is a modified oriented tricolored tree tractable to enumerate. Finally, reducing this bijection to factorizations of a long cycle into two permutations, we get the analogue formula for the corresponding generating series. 
\end{abstract}

%\clearpage
%\tableofcontents

\section{Introduction}

\subsection{Generating series for connection coefficients} \label{sec:conncoeff}

In what follows, we denote by $\la=(\la_1,\la_2,\ldots,\la_k) \vdash n$ an integer partition of $n$ and $\ell(\la)=k$ the length or number of parts of $\la$. We also write $\la=[1^{n_1(\lambda)},2^{n_2(\lambda)},\ldots]$ where $n_i(\lambda)$ is the number of parts $i$ in $\lambda$. 
%Thus, $\la=(\la_1,\ldots, \la_k)$ where $\la_1 \geq \cdots \geq \la_k\geq 1$ and $\sum \la_i=n$. 

Let $\mathfrak{S}_n$ be the symmetric group on $n$ elements, and $\mathcal{C}_{\la}$ be the conjugacy class in $\mathfrak{S}_n$ of permutations with cycle type $\la$, where $\la \vdash n$. Given $\la^{(1)},\la^{(2)},\ldots,\la^{(r)},\mu \vdash n$, let $k_{\la^{(1)},\ldots,\la^{(r)}}^{\mu}$ be the number of ordered factorizations in $\mathfrak{S}_n$ of a fixed permutation $\gamma \in \mathcal{C}_{\mu}$ as a product $\al_1\cdots \al_r$ of $r$ permutations $\al_i \in \mathcal{C}_{\la^{(i)}}$. These numbers are called {\em connection coefficients} of the symmetric group. The problem of computing these coefficients has received significant attention and a good account of its history and references can be found in \cite{GS}. We focus on the cases $k_{\la,\mu}^{n}$ and $k_{\la,\mu,\nu}^{n}$: {i.e.} when $r=2$ and $3$, $\mu=(n)$ and $\gamma$ is the long cycle $\gamma_{n}=(1,2,\ldots, n)$.

In addition, for $\la \vdash n$ we use the monomial symmetric function $m_{\la}({\bf x})$ on indeterminates ${\bf x} = (x_1, x_2,\ldots)$ which is the sum of all different monomials obtained by permuting the variables of $x_1^{\la_1}x_2^{\la_2}\cdots$, and the power symmetric function $p_{\lambda}({\bf x})$, defined multiplicatively as $p_{\la}=p_{\la_1}p_{\la_2}\cdots$ where $p_n({\bf x})=m_n({\bf x})=\sum_{i} x_i^n$. Also, if $\la=[1^{n_1(\la)},2^{n_2(\la)},\ldots]$, let $Aut(\la) = \prod_i n_i(\la)!$.

Our combinatorial results can be stated as follows:
\begin{thm} \label{thm1}
The numbers $k_{\la,\mu,\nu}^n$ of factorizations of the long cycle $\ga_n$ into an ordered product of three permutations of types $\la,\mu$, and $\nu$ respectively satisfy: 

\begin{equation} \label{eq1}
\sum_{\la, \mu, \nu \vdash n} {k_{\la,\mu,\nu}^n} p_{\la}({\bf x})p_{\mu}({\bf y})p_{\nu}({\bf z}) =  \sum_{\la,\mu, \nu \vdash n} \frac{n!^2 M^{(n-1)}_{\ell(\la),\ell(\mu),\ell(\nu)}}{\binom{n-1}{\ell(\la)-1}\binom{n-1}{\ell(\mu)-1}\binom{n-1}{\ell(\nu)-1}}m_\la({\bf x})m_\mu({\bf y})m_{\nu}({\bf z}),
\end{equation}
where: % $p_{\lambda}$ and $m_{\lambda}$ are power sum and monomial symmetric functions, and
$$
M^{(n-1)}_{\ell(\la),\ell(\mu),\ell(\nu)} = \binom{n-1}{\ell(\nu)-1}\sum_{g\geq 0} \binom{n-\ell(\mu)}{\ell(\la)-1-g}\binom{n-\ell(\nu)}{g}\binom{n-1-g}{n-\ell(\mu)}.
$$

\end{thm}

%\subsection{Bijective mapping for bicolored partitioned maps that preserves type}

\begin{cor} \cite{MV} \label{cor1} The numbers $k_{\lambda,\mu}^{n}$ of factorizations of the long cycle $\gamma_n$ into an ordered product of two permutations of cycle types $\lambda$ and $\mu$ respectively satisfy:
\begin{equation} \label{eq2}
\sum_{\la, \mu \vdash n} {k_{\la,\mu}^n} p_{\la}({\bf x})p_{\mu}({\bf y}) =  \sum_{\la, \mu \vdash n} \frac{n(n-\ell(\la))!(n-\ell(\mu))!}{(n+1-\ell(\la)-\ell(\mu))!}m_{\la}({\bf x})m_{\mu}({\bf y}),
\end{equation}
%where $p_{\lambda}$ and $m_{\lambda}$ are power sum and monomial symmetric functions.
\end{cor}

We will see in Section \ref{sect:combref} that the coefficients on the right hand sides of \eqref{eq1} and \eqref{eq2} are non-negative integers.  

\begin{rem}\label{rem0}
Equations \eqref{eq1} and \eqref{eq2} can be obtained algebraically using the irreducible characters of the symmetric group, the Murnaghan-Nakayama rule, and symmetric function identities (see \cite{DMJ}). Here, we derive these equations through a bijection. 
\end{rem}

\subsection{Background} \label{S:Back}

%We denote by $\multiset{X}{k}$, the set of $k$-multisets $M$ of $X$. If $X$ is finite then $\multiset{X}{k}=\multiset{\#X}{k}=\binom{\#X+k-1}{k}$.

 In the setting of the connection coefficients $k_{\la^{(1)},\cdots,\la^{(r)}}^{n}$, we define the {\em genus} ${\sf g}(\la^{(1)},\ldots,\la^{(r)})$ of the partitions $\la^{(i)}$  by the equation 
\begin{equation} \label{eq:genus}
\ell(\la^{(1)})+\cdots + \ell(\la^{(r)})=(r-1)n+1-2{\sf g}(\la^{(1)},\cdots,\la^{(r)}).
\end{equation}
We can take ${\sf g}$ to be a non-negative integer, since otherwise it is easy to show that $k_{\la^{(1)}\cdots\la^{(r)}}^n=0$. 
 
Except for special cases there are no closed formulas for the connection coefficients $k^n_{\la^{(1)},\ldots,\la^{(r)}}$. For instance, using an inductive combinatorial argument B\'{e}dard and Goupil \cite{BG} found a formula for $k_{\la,\mu}^{n}$ in the case ${\sf g}(\la,\mu)=0$. This was extended by Goulden and Jackson \cite{GJ92} to evaluate $k^n_{\la^{(1)},\ldots,\la^{(r)}}$ in the case ${\sf g}(\la^{(1)},\ldots, \la^{(r)})=0$ via a bijection with a set of ordered rooted $r$-cacti on $n$ $r$-gons. Later, using characters of the symmetric group and a combinatorial development, Goupil and Schaeffer \cite{GS} derived an expression for connection coefficients of arbitrary genus as a sum of positive terms (see Biane \cite{PB} for a succinct algebraic derivation; Poulalhon and Schaeffer \cite{PS} and Irving \cite{JI} for further generalizations). As a general rule, these developments are quite intricate and the formulas obtained are rather complicated.
% Recently, Chapuy \cite{GCT} and \cite{GC} found a remarkable combinatorial identity for $\sum_{\la,\mu, \ell(\la)=a,\ell(\mu)=b} k^n_{\la,\mu}$ by identifying $2{\sf g}$ special elements in the cycles, called {\em trisections}, and slicing the cycles at these to obtain lower genus factorizations of the long cycle. 

Interestingly, if we consider the generating series for the coefficients $k_{\la^{(1)},\cdots,\la^{(r)}}^{n}$ as in the LHS of \eqref{eq1}, the coefficients of their expansion in the basis of monomial symmetric functions, as in the RHS of \eqref{eq1}, can be computed in closed form thanks to a result by Jackson \cite{DMJ} obtained algebraically using the theory of the irreducible characters of the symmetric group. There are direct bijections for a variant of the case of two factors (i.e. $r=2$) like the celebrated Harer-Zagier formula \cite{HZ}: see Lass \cite{L}, Goulden and Nica \cite{GN}, and Bernardi \cite{B}.  In this paper we follow this approach and introduce the notion of {\em partitioned tricolored (bicolored) $3$-cacti (maps)} of
given type, refining the work of Schaeffer and Vassilieva in \cite{SV} and Vassilieva in \cite{V}, and use a purely combinatorial argument to derive the explicit generating series for $k^n_{\la,\mu,\nu}$ and $k^n_{\la,\mu}$ in Equations \eqref{eq1} and \eqref{eq2} respectively.
%derive a novel simpler formula for numbers related to  $k^n_{\la,\mu,\nu}$ and $k^n_{\la,\mu}$. 
%Also in the genus zero case, the argument reduces to the bijection in \cite{GJ92} for $k^n_{\la,\mu}$ but differs to the one for $k^n_{\la,\mu,\nu}$. 

\subsection{Outline of paper}
The paper is organized as follows: in Section \ref{sect:combref} we introduce the {\em partitioned $3$-cacti} and the {\em cactus trees} (the enumeration of the latter is postponed to Section  \ref{sect:pfpropnumcact})  and relate them via a bijection $\Theta$ described in Section \ref{sect:descbij}. Finally, in Section \ref{sect:pfcor1} we prove Corollary \ref{cor1}.
%The paper is organized i

%\section{Link with maps and partitioned maps of specified type}
%
%\subsection{Unicellular partitioned bicolored maps}
%
%
%\subsection{Connection between $C(\la,\mu,\nu)$ and $k^n_{\la,\mu,\nu}$}

\medskip

\section{Combinatorial reformulation} \label{sect:combref}
\subsection{Cacti and partitioned cacti} \label{cactisect:combref}
Factorizations in the symmetric group counted by $k_{\lambda, \mu, \nu}^n$ admit a direct interpretation in terms of {\it unicellular $3$-constellations} also named {\it $3$-cacti} with white, black, and grey vertices of respective degree distribution $\lambda$, $\mu$, and  $\nu$. Within a topological point of view, $3$-cacti are specific {\it maps} which in turn are $2-$cell decompositions of an oriented surface into a finite number of vertices ($0-$cells), edges ($1-$cells) and faces ($2-$cells) homeomorphic to open discs (see \cite{LZ} for more details about maps and their applications). Maps are defined up to a homeomorphism of the surface that preserves its orientation, the type of cells, and incidences in the graph. $3$-cacti are maps with black faces and one white face (thus the term {\it unicellular}) such that: (i) each edge separates a black face and the white face and (ii) all the black faces are triangles each composed of exactly a white, a black, and a grey vertex following each other in clockwise order. As a consequence, the degree of the white face is a multiple of $3$. Often, cacti refer to planar maps (embedded in an orientable surface of genus $0$). In this paper we assume that they can be embedded in an orientable surface of any genus. Besides, we consider only {\it rooted} cacti, i.e. cacti with a marked black face. We assume as well that each black triangle is labeled with an index in $\{1,2,\ldots,n\}$ with the convention that the marked triangle is labeled $1$. In what follows, we define the {\it degree} of a vertex in a cactus as the number of triangles it belongs to, and the degree distribution of the vertices of a given color is the integer partition of $n$ formed by the degrees of all the vertices of this color. 

The next classical result (see \cite{LZ}) relates rooted $3$-cacti with factorizations of the long cycle $\gamma_n=(1,2,\ldots,n)$.

\begin{prop}
\label{propo}
Rooted $3$-cacti with $n$ black triangles are in bijection with $3$-tuples $(\alpha_1$, $\alpha_2$, $\alpha_3)$  of permutations in $\mathfrak{S}_n$ such that $\alpha_1\alpha_2\alpha_3 = \gamma_n$.  Under this bijection the white (black and grey, resp.) vertices correspond to cycles of $\pi_1$ ($\pi_2$ and $\pi_3$, resp.).
\end{prop}

A sketch of the proof  of this classical result can be found in \cite{V}. Each white, black, or grey vertex of a given $3$-cacti corresponds to a cycle of permutation $\al_1$, $\al_2$, or $\al_3$ respectively, and the cycle is encoded by the local {counter-clockwise} order of the triangles around the vertex. The fact that $\al_1\al_2\al_2=\gamma_n$ corresponds to saying that traversing the map starting on the white vertex of the triangle labeled $1$ and keeping the white face on the {\em right} we visit, in order, the white-black edges belonging to triangles labeled $1,2,\ldots,n$.  A consequence of Proposition~\ref{propo} is that for integer partitions $\lambda,\mu,\nu$ of $n$,  the number of $3$-cacti of degree distribution $\la$, $\mu$, $\nu$ is  the number $k_{\lambda,\mu,\nu}^n$ of factorizations defined in Section~\ref{sec:conncoeff}. Moreover, by the Euler-characteristic, the {genus} of the underlying surface of the $3$-cacti of degree distribution $\la, \mu, \nu$ is given by Equation \ref{eq:genus}. Figure 1 %{fig:cacdonut}
shows a $3$-cactus embedded on the sphere (genus $0$) and a $3$-cactus embedded on a torus (genus $1$).

\begin{figure}[h]
\begin{center}
\raisebox{.2in}{\label{fig:cacdonut}
\includegraphics[width=50mm]{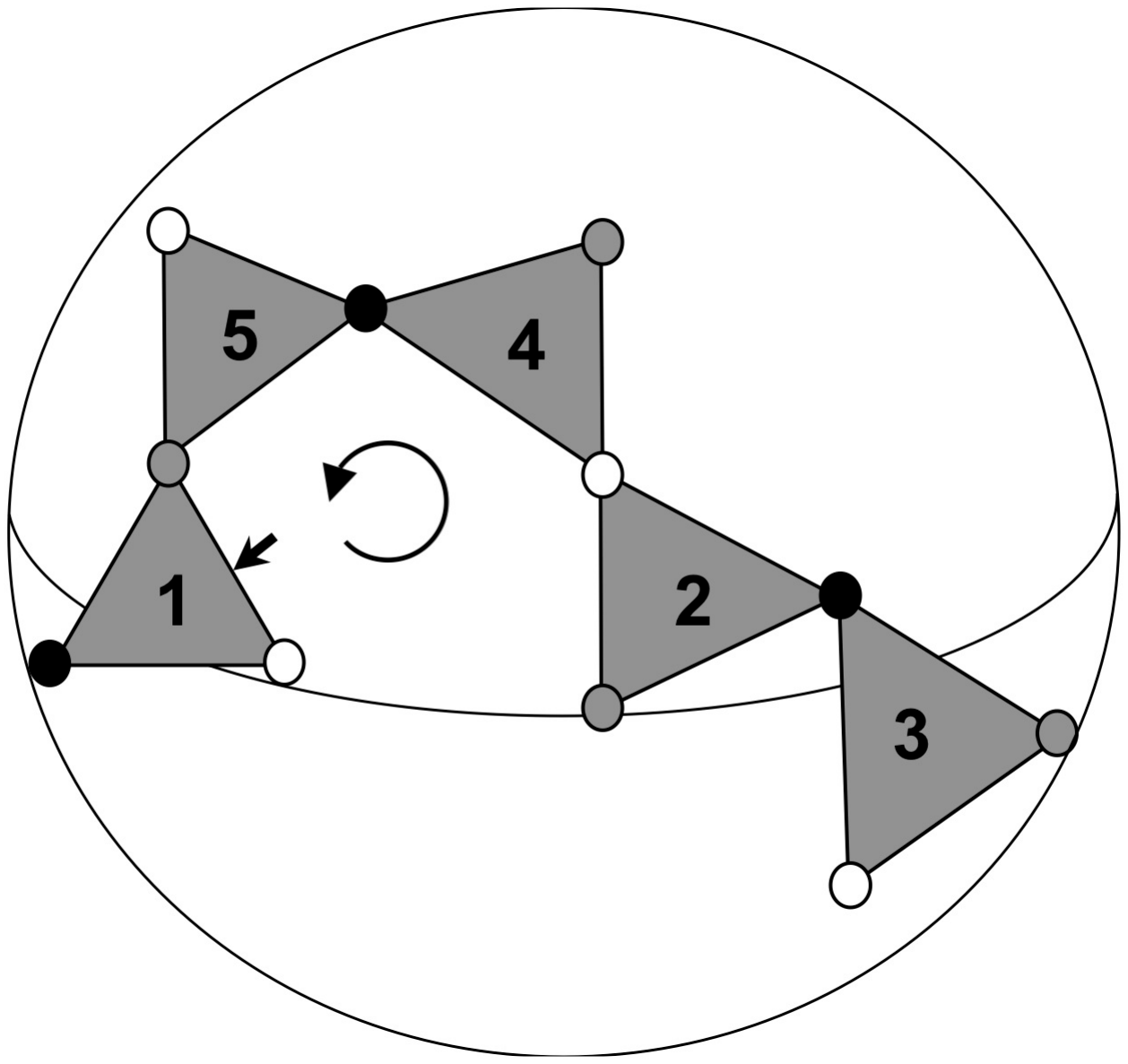}
}
\includegraphics[width=60mm]{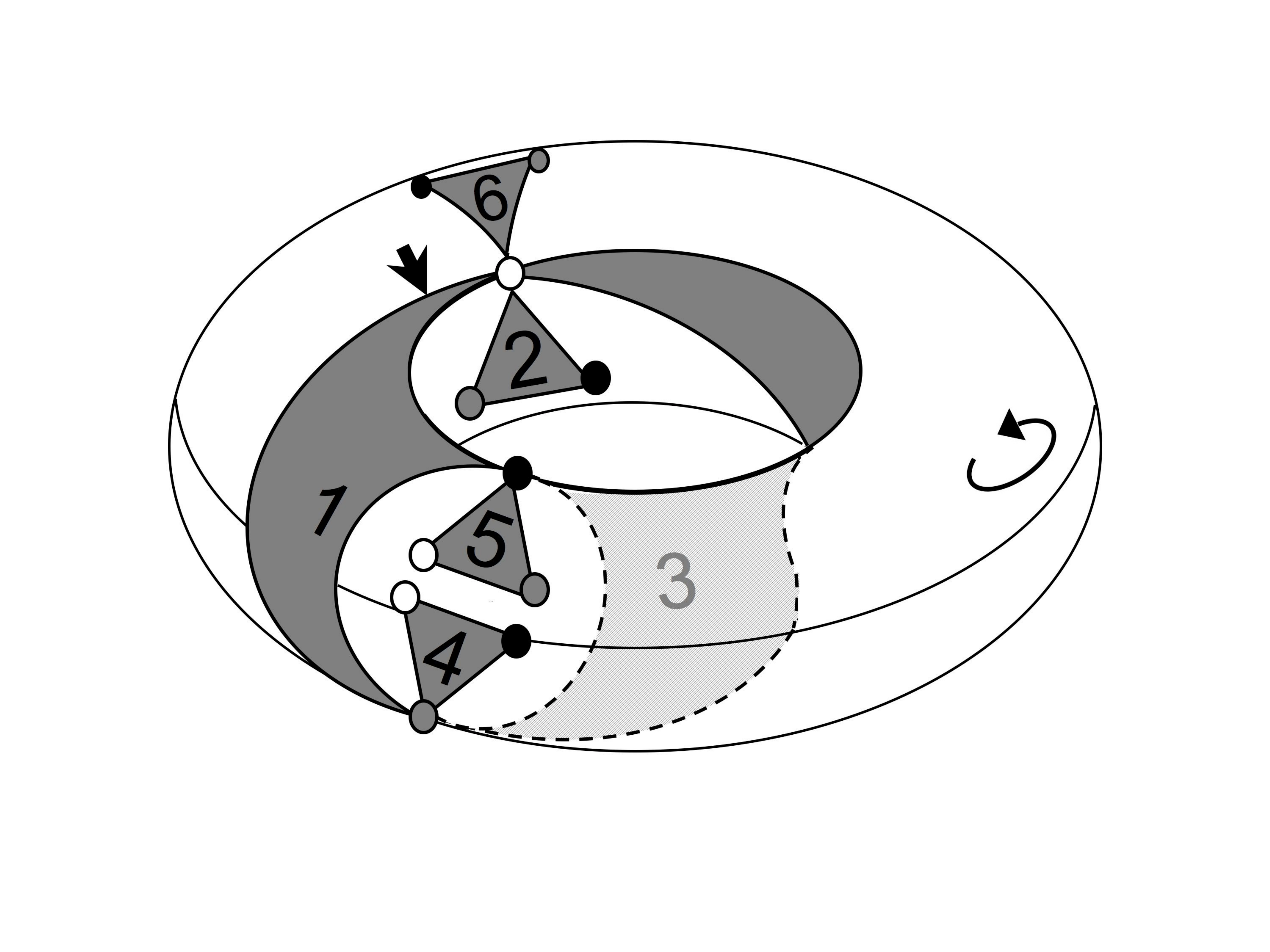}
\caption{Two examples of rooted $3$-cacti embedded in a surface of genus 0 (left) and genus 1 (right)}
\end{center} 

\end{figure}

\begin{exm}

The cactus on the left hand side in Figure 1 %\ref{fig:cacdonut}
can be associated to the three permutations:
\begin{eqnarray}
\alpha_1 &=& (1)(24)(3)(5)\\
\alpha_2 &=& (1)(23)(45)\\
\alpha_3 &=& (15)(2)(3)(4)
\end{eqnarray}
The degree distribution of this cactus is $\la = [1^3,2^1]$, $\mu = [1^1,2^2]$, $\nu = [1^3,2^1]$.\\
The cactus on the right hand side corresponds to the three permutations:
\begin{eqnarray}
\alpha_1 &=& (1236)(4)(5)\\
\alpha_2 &=& (153)(2)(4)(6)\\
\alpha_3 &=& (134)(2)(5)(6)
\end{eqnarray}
The degree distribution of this cactus is $\la = [1^2,4^1]$, $\mu = [1^3,3^1]$, $\nu = [1^3,3^1]$.
\end{exm}

\noindent {\bf{Partitioned 3-cacti}}\\
Cacti with a given vertex degree distribution are in general non-planar and non-recursive objects, and no direct bijection is known to compute their cardinality.
%As non recursive and non planar objects, the cardinality of general sets of cacti with a given vertex degree distribution cannot be computed directly. 
However, we provide an interpretation of the formal power series $\sum_{\la, \mu, \nu \vdash n} {k_{\la,\mu,\nu}^n} p_{\la}({\bf x})p_{\mu}({\bf y})p_{\nu}({\bf z})$ as a generating function for {\it partitioned cacti}. We are able to give an explicit formula for this generating function by introducing a new bijective mapping for partitioned cacti. Intuitively, partitioned cacti are rooted $3$-cacti where the set of vertices of each color are partitioned into blocks. Such objects have also been widely studied, for instance by Lass \cite{L} and Bernardi \cite{B} under the term of colored maps; and by Goulden and Nica \cite{GN}, Schaeffer and Vassilieva in \cite{SV}, and Vassilieva in \cite{V} as partitioned maps or cacti.
%Partitioned maps and cacti have also been studied . 

To define partitioned cacti we use $\pi$ to denote a set partition of a set of $n$ elements with blocks $\{\pi^1,\ldots, \pi^p\}$. The type of a set partition, $\ty(\pi) \vdash n$, is the integer partition of $n$ obtained by considering the cardinalities of the blocks of $\pi$. We are now ready for the definition:
\begin{defn}(Partitioned 3-cactus)
A {\it partitioned $3$-cactus} is a $4$-tuple $(\pi_{1}, \pi_{2}, \pi_{3},\kappa)$ such that $\kappa$ is a rooted $3$-cactus with $n$ triangles, and $\pi_1$, $\pi_2$ and $\pi_3$ are set partitions on the set of white, black, and grey vertices respectively. By abuse of notation, hereinafter we view $\pi_1$, $\pi_2$, and $\pi_3$ as set partitions on $\{1,2,\ldots,n\}$ where a block is composed of the labels of the triangles incident to the vertices contained in the block.  In what follows, the blocks of $\pi_{1}, \pi_{2}$, and $\pi_{3}$ are denoted $\pi_{1}^{(i)}, \pi_{2}^{(j)}$, and $\pi_{3}^{(k)}$ with the only restriction that $1 \in \pi_{1}^{(\ell(\la))}$.

For $\lambda,\mu, \nu \vdash n$, we let $\mathcal{C}(\lambda,\mu,\nu)$ be the set of partitioned cacti $(\pi_{1}, \pi_{2}, \pi_{3},\kappa)$ where the set partitions $\pi_1$, $\pi_2$, and $\pi_3$ of $\{1,2,\ldots,n\}$ have type $\lambda$, $\mu$, and $\nu$ respectively. Let $C(\la, \mu, \nu)= |\mathcal{C}(\la, \mu, \nu)|$.  
%For $\la$, $\mu$, $\nu \vdash n$, we let $\mathcal{C}(\la, \mu, \nu)$ be the set of partitioned $3$-cacti $(\pi_{1}, \pi_{2}, \pi_{3},\kappa)$ such that $\kappa$ is a $3$-cactus with $n$ triangles, $\pi_1$ (resp. $\pi_2$, $\pi_3$) a set partition on the set of white (black and grey, resp.) vertices. By abuse of notation, $\pi_1$, $\pi_2$ and $\pi_3$ are viewed as set partitions on $\{1,2,\ldots,n\}$ of type $\lambda$ (resp. $\mu,\nu$) where a block is composed of the labels of the triangles adjacent to the vertices contained in the block.  In what follows, the blocks of $\pi_{1}, \pi_{2}$, and $\pi_{3}$ are denoted $\pi_{1}^{(i)}, \pi_{2}^{(j)}$, and $\pi_{3}^{(k)}$ with the only restriction that $1 \in \pi_{1}^{(\ell(\la))}$. Let $C(\la, \mu, \nu)= |\mathcal{C}(\la, \mu, \nu)|$.  
\end{defn}

%Note that in terms of factorizations, partitioned cacti can be viewed as ordered factorizations of $\gamma_n$ where we label the cycles of permutations, using distinct sets of labels for each permutation, allowing repeated labels within the permutation. Then the parts of $\lambda$ are the number of elements in $[n]$ belonging to cycles of $\al_1$ with the same label. An equivalent statement holds for $\mu$ and $\nu$.

\begin{rem}
Following Proposition \ref{propo}, we can interpret the cactus $\kappa$ as a $3$-tuple of permutations $(\alpha_1,\alpha_2,\alpha_3)$ where $\alpha_3=\alpha_{2}^{-1}\circ\alpha_{1}^{-1}\circ\gamma_{n}$. As a result, partitioned cacti in $\mathcal{C}(\la, \mu, \nu)$ can be seen as 5-tuples $(\pi_{1}, \pi_{2}, \pi_{3}, \alpha_{1}, \alpha_{2})$ where $\pi_{1}$,$\pi_{2}$, and $\pi_{3}$
are set partitions of $\{1,2,\ldots,n\}$ of types $\lambda$, $\mu$, and $\nu$ respectively with the property that: for $k=1,2,3$, if an integer $l$ of a given cycle $c$ of $\al_k$ belongs to a given block of $\pi_k$ then all the integers in the cycle $c$ also belong to this block.
%composed of set partitions $\pi_{1}$,$\pi_{2}$, and $\pi_{3}$
%with respective type $\la$, $\mu$, and $\nu$; and permutations $\al_1,\al_2 \in \mathfrak{S}_{n}$ such that:
%\begin{compactitem}
%\item {\tt } each block of $\pi_{1}$ is the union of cycles of $\alpha_{1}$,
%\item {\tt } each block of $\pi_{2}$ is the union of cycles of $\alpha_{2}$,
%\item {\tt } each block of $\pi_{3}$ is the union of cycles of $\alpha_{3}=\alpha_{2}^{-1}\circ\alpha_{1}^{-1}\circ\gamma_{n}$.
%%where $\gamma_{N}• = (1\, 2\, \ldots\, N)$.
%\end {compactitem}
\end{rem}

%
%\\
%{\bf{Graphical Representation}}\\
%Graphically, a  $5$-tuple  $(\pi_{1}, \pi_{2}, \pi_{3}, \alpha_{1}, \alpha_{2})$ corresponds to
%a $3$-cactus with N triangles:
%\begin{itemize}
%\item {\tt } the cycles of $\alpha_{1}$ describe the white vertices of the cactus,
%\item {\tt } the cycles of $\alpha _{2}$ describe the black vertices,
%\item {\tt } the cycles of $\alpha_{3} = \alpha_{2} ^{-1}\circ \alpha_{1} ^{-1}\circ \gamma_{N}$ describe the grey vertices,
%\item{\tt }  $\pi_{1}$ partitions the white vertices into $p_{1}$ subsets,
%\item{\tt }  $\pi_{2}$ partitions the black vertices into $p_{2}$ subsets ,
%\item{\tt }  $\pi_{3}$ partitions the grey vertices into $p_{3}$ subsets.
%\end {itemize}

\begin{exm}
\label{ex: pc1}
 Take the cactus on the right hand side of Figure 1 %\ref{fig:cacdonut}
 and add partitions $\pi_{1}=\{ \pi_{1}^{(1)}, \pi_{1}^{(2)}\}, \pi_{2}=\{\pi_{2}^{(1)}, \pi_{2}^{(2)}\}$,
$\pi_{3}=\{\pi_{3}^{(1)}, \pi_{3}^{(2)}, \pi_{3}^{(3)}\}$, where
$\pi_{1}^{(1)}=\{4, 5\},\,\,\,\pi_{1}^{(2)}=\{1, 2, 3, 6\}; \; \pi_{2}^{(1)}=\{1, 3, 4, 5\},\,\,\, \pi_{2}^{(2)}=\{2, 6\};\; \pi_{3}^{(1)}=\{1, 3, 4, 6\},\,\,\, \pi_{3}^{(2)}=\{2\},\,\,\, \pi_{3}^{(3)}=\{5\}$. This gives a partitioned cactus $(\pi_{1}, \pi_{2}, \pi_{3}, \alpha_{1}, \alpha_{2})  \in \mathcal{C}({[2^1,4^1],[2^1,4^1], [1^2,4^1]})$ depicted in Figure \ref{fig:ex}.
Similarly to \cite{SV}, we associate a particular shape to each of the blocks of the partitions.

\begin{figure}[h]
\begin{center}
\begin{minipage}[b]{0.2\linewidth}
\includegraphics[height=3.2cm]{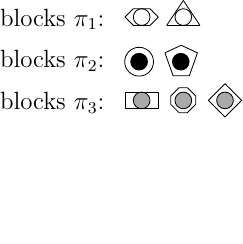}
\end{minipage}
\begin{minipage}[b]{0.7\linewidth}
\includegraphics[width=75mm]{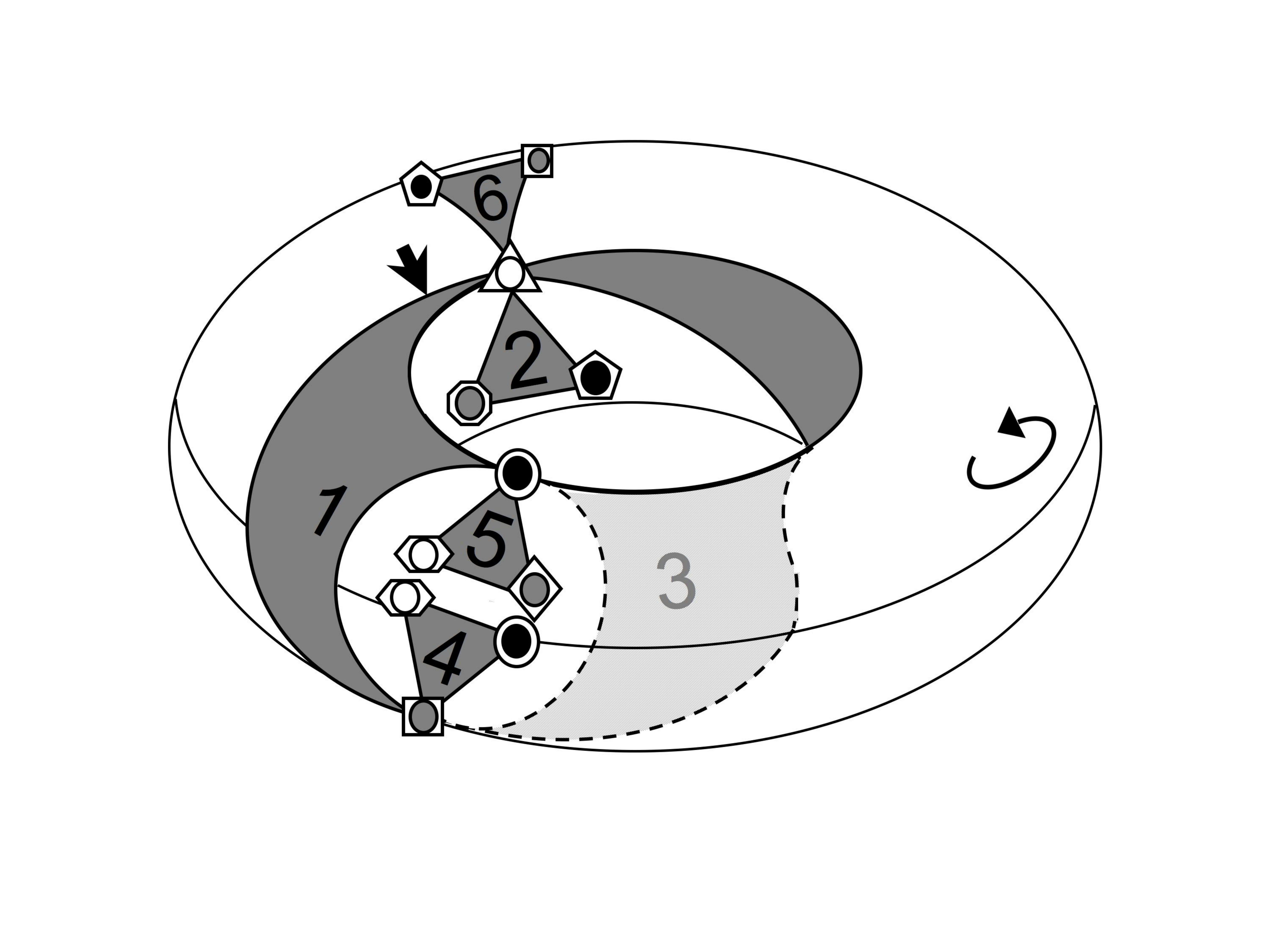}
\end{minipage}
\caption{Illustration of the Partitioned $3$-Cactus from Example~\ref{ex: pc1}.}
\label{fig:ex}
\end{center} 
\end{figure}

\end{exm}

\noindent {\bf{Link between cacti and partitioned cacti}}\\
Consider the partial order on integer partitions given by refinement. That is $\la \preceq \mu$ if and only if the parts of $\mu$ are unions of parts of $\la$, and we say that $\mu$ is {\em coarser} than $\la$. If $\la \preceq \mu$ let $\ov{R}_{\la,\mu}$ be the number of ways to {\em coarse} $\la$  to obtain $\mu$. For example, if $\la =[1
^2,2^2]$ and $\mu=[1,2,3]$ then $\ov{R}_{\la,\mu}=4$.
% since in $\la$ any of the two $1$-blocks can merge with any of the two $2$-blocks.
It is well known that 
$p_{\la} = \sum_{\mu \succeq \la} Aut(\mu)\ov{R}_{\la,\mu} m_{\mu}$ \cite[Prop. 7.7.1]{EC2}.

We use this partial order on integer partitions to obtain an immediate relation between $C(\la,\mu,\nu)$ and $k_{\la,\mu,\nu}^{n}$.
 
\begin{prop} \label{prop1} For partitions $\rho,\delta,\epsilon \vdash n$ we have :

\begin{equation} \label{eq:prop1}
C(\rho,\delta,\epsilon) = \sum_{\la \preceq \rho, \mu \preceq \delta, \nu \preceq \epsilon} \ov{R}_{\la\rho}\ov{R}_{\mu\delta}\ov{R}_{\nu\epsilon}k_{\la,\mu,\nu}^n. 
\end{equation}

\end{prop}

\begin{proof}
Let $(\pi_1, \pi_2,\pi_3,\al_1,\al_2) \in \mathcal{C}(\rho,\delta,\epsilon)$. If $\al_1 \in \mathcal{C}_{\la}$, $\al_2 \in \mathcal{C}_{\mu}$, and $\al_3 = \al_2^{-1}\al_1^{-1}\ga_n \in \mathcal{C}_{\nu}$  then by the definition of the partitioned cacti, we have that $\ty(\pi_1)=\rho \succeq \la$, $\ty(\pi_2)=\delta \succeq \mu$, and $\ty(\pi_3)=\epsilon \succeq \nu$. Thus, if we further refine $C(\rho,\delta,\epsilon)$ by the cycle types of the permutations, {i.e.} if
$$
\mathcal{C}_{\la,\mu,\nu}(\rho,\delta,\epsilon) = \{ (\pi_1,\pi_2,\pi_3,\al_1,\al_2) \in \mathcal{C}(\rho,\delta,\epsilon) ~|~ (\al_1,\al_2,\al_2^{-1}\al_1^{-1}\ga_n) \in \mathcal{C}_{\la} \times \mathcal{C}_{\mu} \times \mathcal{C}_{\nu}\},
$$
then $\mathcal{C}(\rho,\delta,\epsilon) = \bigcup_{\la \preceq \rho, \mu \preceq \delta, \nu \preceq \epsilon} \mathcal{C}_{\la,\mu,\nu}(\rho,\delta,\epsilon) $ where the union is disjoint. Finally, if $$C_{\la,\mu,\nu}(\rho,\delta,\epsilon) = |\mathcal{C}_{\la,\mu,\nu}(\rho,\delta,\epsilon)|$$
then it is easy to see that $C_{\la,\mu,\nu}(\rho,\delta,\epsilon)=  \ov{R}_{\la\rho}\ov{R}_{\mu\delta}\ov{R}_{\nu\epsilon}k_{\la,\mu,\nu}^n$. 
\end{proof}

\noindent Using $p_{\la} = \sum_{\mu \succeq \la} Aut(\mu)\ov{R}_{\la,\mu} m_{\mu}$ Proposition \ref{prop1} is equivalent to:
% in terms of the generating function $\psi_n(x,y)$. 

\begin{equation} \label{ngs}
\sum_{\la, \mu,\nu \vdash n} k_{\la,\mu,\nu}^n p_{\la}({\bf x})p_{\mu}({\bf y})p_{\nu}({\bf z}) = \sum_{\la,\mu,\nu \vdash n} Aut(\la)Aut(\mu)Aut(\nu)\,C(\la,\mu,\nu)\, m_{\la}({\bf x})m_{\mu}({\bf y})m_{\nu}({\bf z}) 
\end{equation}

In the special case when we have partitions $\rho, \delta$ and $\epsilon$ of $n$ where $\ell(\rho) + \ell(\delta) + \ell(\epsilon)=2n+1$, the following proposition holds:

\begin{prop}[\cite{GJ92}] \label{proprem2}
For partitions $\rho, \delta$ and $\epsilon$ of $n$ where $\ell(\rho) + \ell(\delta) + \ell(\epsilon)=2n+1$ we have that $C(\rho,\delta,\epsilon)=k^n_{\rho,\delta,\epsilon}=n^2(\ell(\rho)-1)!(\ell(\delta)-1)!(\ell(\epsilon)-1)!/{Aut(\rho)Aut(\delta)Aut(\epsilon)}$.
\end{prop} 

\begin{proof}
Let $(\pi_1,\pi_2,\pi_3,\al_1,\al_2) \in \mathcal{C}(\rho,\delta,\epsilon)$ with $\ell(\rho)+ \ell(\delta)+\ell(\epsilon)=2n+1$, and $\al_3=\al_2^{-1}\al_1^{-1}\gamma_n$. If $\al_1\in \mathcal{C}_{\la}$, $\al_2\in \mathcal{C}_{\mu}$, $\al_3\in \mathcal{C}_{\nu}$, then $\ell(\la) + \ell(\mu)+\ell(\nu)=2n+1-2{\sf g}(\la,\mu,\nu)\leq 2n+1$. But $\ell(\la) \geq \ell(\rho)$, $\ell(\mu)\geq \ell(\delta)$, and $\ell(\nu)\geq \ell(\epsilon)$ therefore $\rho=\la$, $\delta=\mu$, and $\epsilon=\nu$; and $\pi_1$, $\pi_2$, and $\pi_3$ are the underlying set partitions in the cycle decompositions of $\al_1, \al_2$, and $\al_3$ respectively. Thus $C(\delta,\rho,\epsilon)=k_{\rho,\delta,\epsilon}^{n}$. But, as shown in  \cite[Thm. 2.2]{GJ92}, $k^n_{\rho,\delta,\epsilon}=n^2(\ell(\rho)-1)!(\ell(\delta)-1)!(\ell(\epsilon)-1)!/{Aut(\rho)Aut(\delta)Aut(\epsilon)}$ since the genus ${\sf g}(\rho,\delta,\epsilon)=0$. As a result, $C(\rho,\delta,\epsilon) = n^2(\ell(\rho)-1)!(\ell(\delta)-1)!(\ell(\epsilon)-1)!/{Aut(\rho)Aut(\delta)Aut(\epsilon)}$
\end{proof}

As mentioned before, explicit computation of the right hand side of equation \eqref{ngs} is made possible thanks to a new bijective description of partitioned cacti of given type which is a refinement of a bijection in \cite{SV} and \cite{V}. Partitioned cacti are indeed in one-to-one correspondence with particular sets of cactus trees (and three additional simple combinatorial objects), which are recursive planar objects whose number one can compute with classical methods like Lagrange inversion. Next we define such trees, in Section~\ref{sect:pfpropnumcact} we compute the number of such trees.

\subsection{Cactus trees} \label{subsec:ct}

\begin{figure}[h]
\begin{center}
\includegraphics[width=40mm]{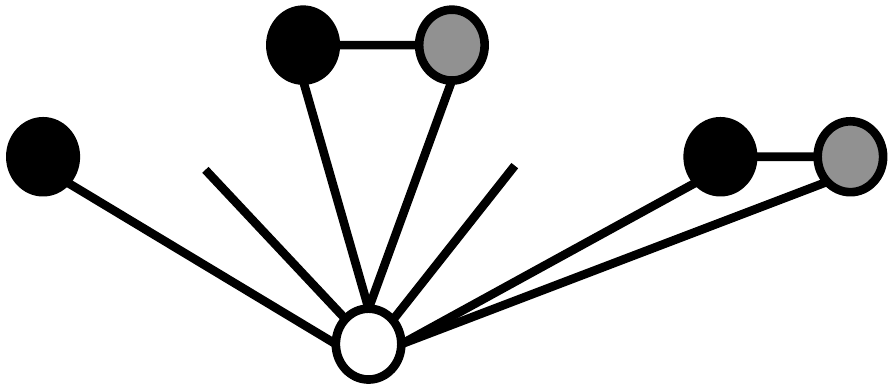}
\caption{Example of three types of possible children, the ordered sequence of children attached to the white vertex is: edge, thorn, triangle, thorn, triangle}
\label{children}
\end{center} 
\end{figure}

Before we state the actual definition of the tree structure used as the main ingredient in the proof of Theorem \ref{thm1}, we give preparatory explanations.
Ordered trees are non cyclic graphs usually defined recursively as a root vertex $v$ and an ordered sequence (possibly the empty set) of ordered trees, called descending trees, each having its root vertex connected to $v$ by an edge. The root of a descending vertex is called a descending vertex. The root vertices of the descending trees of a given vertex are considered as its children. Although it follows the same kind of recursive definition, the tree structure we introduce has the following differences:
\begin{itemize}
\item vertices are of three different colors, say white, black and grey;
\item the ordered sequence of children of a given vertex is not composed of vertices connected to it through an edge. A child can be: (i) a {\em thorn} or half edge, i.e an edge, connecting this given vertex to no other (as a result, no descending tree is attached to this kind of child), (ii) a {\em full edge}, i.e. an edge connecting the given vertex to a descending one with the restriction that only a black (resp. grey, white) vertex can be connected this way to a white (resp. black, grey) one, (iii) a {\em triangle} connecting the given vertex to two descending ones with the restriction that only a black and grey (resp. grey and white, white and black) can be connected this way to a white (resp. black, grey) one. Triangles are made of three edges connecting the two descending vertices to the ascending one and the two descending vertices between themselves. The two descending vertices are the roots of two descending trees. The three kinds of children are illustrated on Figure \ref{children}.
\end{itemize}

We are now ready to state the definition of the cactus trees:
\begin{defn}[Cactus Tree] 
Let $\widetilde{\mathcal{CT}}(p_1,p_2,p_3,g,w,b)$ be the set of {\em cactus trees} $\widetilde{\tau}$ with $p_1$ white vertices, $p_2$ black vertices, and $p_3$ grey vertices, $g$ triangles children of grey vertices, $w$ triangles children of white vertices,  and $b$ triangles children of black vertices such that:
\begin{itemize}
\item[(i)] the root of $\widetilde{\tau}$ is a white vertex, 
\item[(ii)] the ordered set of children of each white (resp. black, grey) vertex consists of three kinds of objects: thorns; full edges connecting this white (resp. black, grey) vertex to a black (resp. grey, white) one; triangles connecting this white (resp. black, grey) vertex to both a black and a grey (resp. grey and white, white and black) one, 
\item[(iii)] the edge connecting a white (resp. black, grey) vertex to the black (resp. grey, white) one in a descending triangle comes before the one connecting it to the grey (resp. white, black) vertex according to the ordering of the children of this white (resp. black, grey) vertex.
\end{itemize}

%consisting of an underlying cactus tree $\tau \in \mathcal{CT}(p_1,p_2,p_3,g,w,b)$, and in addition it has: $n+1-p_1-p_2+g$ thorns connected to white vertices, $n-p_2-p_3+w$ thorns connected to black vertices, and $n+1-p_1-p_3+b$ {\em thorns} connected to grey vertices. The only restriction for the thorns is: \\ (vi) thorns {\em cannot be inside a triangle} of the cactus.
\end{defn}

Within a cactus tree, the degree of a vertex $\mathsf{v}$ is defined by:   
\begin{equation}\label{degdef}
deg(\mathsf{v})=c+\varepsilon,
\end{equation}
where $c$ is the number of children (that can be either thorns, edges or triangles) and $\varepsilon$ is $1$ for a non-root vertex, $0$ otherwise. With this definition of degree, we write the set $\widetilde{\mathcal{CT}}(p_1,p_2,p_3,g,w,b)$ as the disjoint union
$$\widetilde{\mathcal{CT}}(p_1,p_2,p_3,g,w,b)=\bigcup_{\ell(\la)=p_1,\ell(\mu)=p_2,\ell(\nu)=p_3} \widetilde{\mathcal{CT}}(\la,\mu,\nu,g,w,b).$$
 where $\widetilde{\mathcal{CT}}(\la,\mu,\nu,g,w,b)$ is the set of cactus trees in $\widetilde{\mathcal{CT}}(p_1,p_2,p_3,g,w,b)$ with degree distribution $\la,\mu,\nu \vdash n$ of the white, black, and grey vertices respectively. 

In addition, in what follows we will denote by $\mathcal{CT}(p_1,p_2,p_3,g,w,b)$ the set of cactus trees similar to those in $\widetilde{\mathcal{CT}}(p_1,p_2,p_3,g,w,b)$ but without thorns. We define $\mathcal{CT}(\la,\mu,\nu,g,w,b)$ similarly. Moreover, we will use the expression {\em tricolored tree} when only full edges are allowed in the set of children of each vertex. Finally, we may use the integers $(1,2,\ldots,p_1)$ (resp. $(1,2,\ldots,p_2)$, $(1,2,\ldots,p_3)$) to label the white (resp. black, grey) vertices of a given cactus tree or tricolored tree. The resulting object is called labeled cactus tree or labeled tricolored tree respectively. 

\begin{exm}
\label{ex:cactree}
The cactus tree in Figure \ref{tct} belongs to $$\widetilde{\mathcal{CT}}([1^2,2^1,4^1],[1^1,3^1,4^1],[1^1,2^2,3^1],1,1,1)$$
\begin{figure}[h]
\begin{center}
\includegraphics[width=40mm]{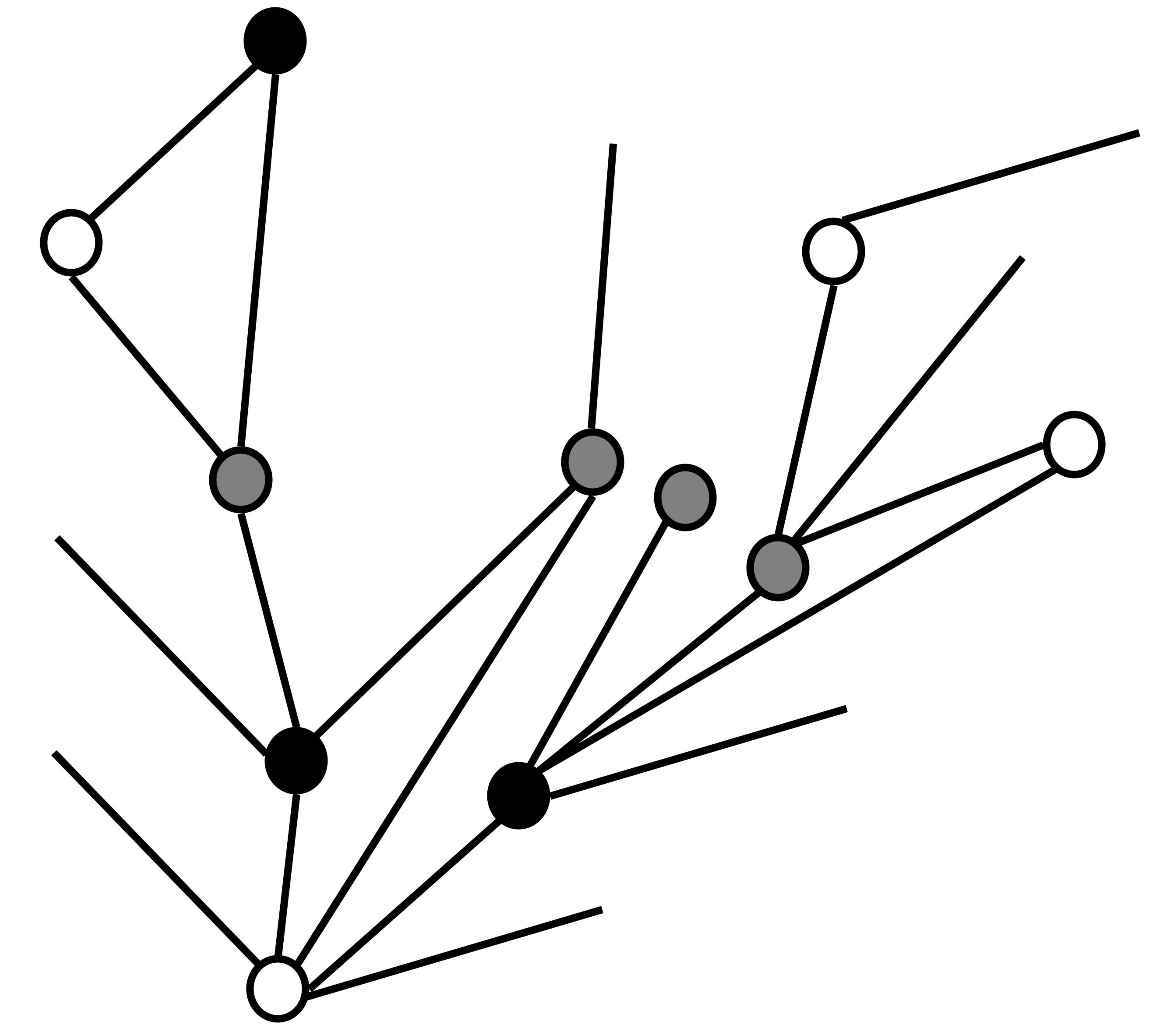}
\caption{Example of a Cactus Tree. The white vertex in the bottom is the root of the tree. The tree has $p_1=3$ white vertices, $p_2=3$ black vertices and $p_3=4$ grey vertices, and three triangles each children of white, black, and grey vertices respectively.}
\label{tct}
\end{center} 
\end{figure}
\end{exm}

\begin{prop}\label{numcact}
The number of cactus trees is:
\begin{multline}
|\widetilde{\mathcal{CT}}(\lambda,\mu,\nu,g,w,b)| = n\cdot \frac{(\ell(\lambda)-1)!(\ell(\mu)-1)!(\ell(\nu)-1)!}{Aut(\lambda)Aut(\mu)Aut(\nu)} \frac{(g(w-\ell(\nu))+\ell(\mu)\ell(\nu))}{(n+1-\ell(\lambda)-\ell(\mu)+g)} \times \\
\times \binom{n-\ell(\lambda)}{w,\ell(\mu)-g-w}\binom{n-\ell(\mu)}{b,\ell(\nu)-w-b}\binom{n-\ell(\nu)}{g,\ell(\lambda)-1-g-b}.
\end{multline}
\end{prop}

The proof of this proposition is carried out using the Lagrange inversion theorem (see e.g. \cite[1.2.13]{GJCE}) and it is postponed to Section \ref{sect:pfpropnumcact}.

\subsection{Reformulation of the main theorem}

Let $\mathcal{OP}_{r}^{(m)}$ be the set of all ordered $r$-subsets of $[m]$. By definition $|\mathcal{OP}_{r}^{(m)}| = (m)_r = m(m-1)\cdots(m-r+1)$. We have the following proposition:

\begin{prop}
\label{MCT}
Theorem \ref{thm1} is equivalent to:
\footnotesize
\begin{equation} \label{eqmct}
C(\la,\mu,\nu) = \sum_{g,w,b \geq 0}|\widetilde{\mathcal{CT}}(\lambda,\mu,\nu,g,w,b)|\,\,|\mathfrak{S}_{n+1-\ell(\la)-\ell(\nu)+b}|\,\,|\mathfrak{S}_{n-\ell(\mu)-\ell(\nu)+w}|\,\,|\mathcal{OP}_{\ell(\nu)-w-b}^{(n+1-\ell(\la)-\ell(\mu)+g)}|
\end{equation}
\normalsize
\end{prop}

\begin{proof}
According to Equation \ref{ngs}, Theorem \ref{thm1} is equivalent to the equality $$C(\la,\mu,\nu) = \frac{n!^2}{Aut(\la)Aut(\mu)Aut(\nu)}\frac{1}{\binom{n-1}{\ell(\la)-1}\binom{n-1}{\ell(\mu)-1}\binom{n-1}{\ell(\nu)-1}}{M}^{(n-1)}_{\ell(\la),\ell(\mu),\ell(\nu)}$$ After basic simplifications on the binomial coefficients, a {\em summand} in the RHS of Equation \eqref{eqmct} reduces to:

\footnotesize
\begin{eqnarray}
\nonumber \frac{(n-1)!^2(\ell(\mu)\ell(\nu)+g(w-\ell(\nu)))}{\binom{n-1}{\ell(\la)-1}\binom{n-1}{\ell(\mu)-1}\binom{n-1}{\ell(\nu)-1}Aut(\lambda)Aut(\mu)Aut(\nu)} \binom{n}{w,g,b,\ell(\la)-1-g-b,\ell(\mu)-g-b,\ell(\nu)-w-b}
\end{eqnarray}
\normalsize
Then we sum over $g$, $w$, and $b$ the terms depending on these parameters. Arranging properly the terms depending on $w$ and $b$, and simplifying sums on these two parameters thanks to the Vandermonde's convolution formula, we obtain:
\begin{multline*}
 \sum_{g,w,b}(\ell(\mu)\ell(\nu)+g(w-\ell(\nu)))\binom{n}{w,g,b,\ell(\la)-1-g-b,\ell(\mu)-g-b,\ell(\nu)-w-b}=\\
= n^2\binom{n-1}{\ell(\nu)-1}\sum_{g}\binom{n-\ell(\mu)}{\ell(\la)-1-g}\binom{n-\ell(\nu)}{g}\binom{n-1-g}{n-\ell(\mu)}.
\end{multline*}
Which leads directly to the desired result.
\end{proof}

As a direct consequence of Proposition \ref{MCT}, Theorem \ref{thm1} reduces to:
\begin{thm}
\label{thm:bij}
There is an explicit bijection $\Theta^n_{\la,\mu,\nu}$ between partitioned 3-cacti in  $\mathcal{C}(\la,\mu,\nu)$ and tuples $(\widetilde{\tau},\sigma_1,\sigma_2,\chi)$, where 
\begin{align*}
&\widetilde{\tau} \in \widetilde{\mathcal{CT}}(\la,\mu,\nu,g,w,b),\\
 &\sigma_1 \in \mathfrak{S}_{n+1-\ell(\la)-\ell(\nu)+b},\\
 &\sigma_2 \in \mathfrak{S}_{n-\ell(\mu)-\ell(\nu)+w},\\
 &\chi \in \mathcal{OP}_{\ell(\nu)-w-b}^{(n+1-\ell(\la)-\ell(\mu)+g)}
 \end{align*}
for some $g,w,b\geq 0$.
\end{thm}

The next section is devoted to proving this theorem by describing the bijection $\Theta^n_{\la,\mu,\nu}$.

\section{Description of the bijection} \label{sect:descbij}

\subsection{Additional definitions}
Before we get to the description of $\Theta^n_{\la,\mu,\nu}$, we need two additional ingredients: a linear order on the white (black and grey, resp.) vertices and their children (as defined in the beginning of Section~\ref{subsec:ct}) which we call white {\em reverse level traversal (RLT)} (black and grey reverse level traversal, resp.), and partial permutations. 

\begin{figure}[h]
\subfigure[]{
\includegraphics[width=50mm]{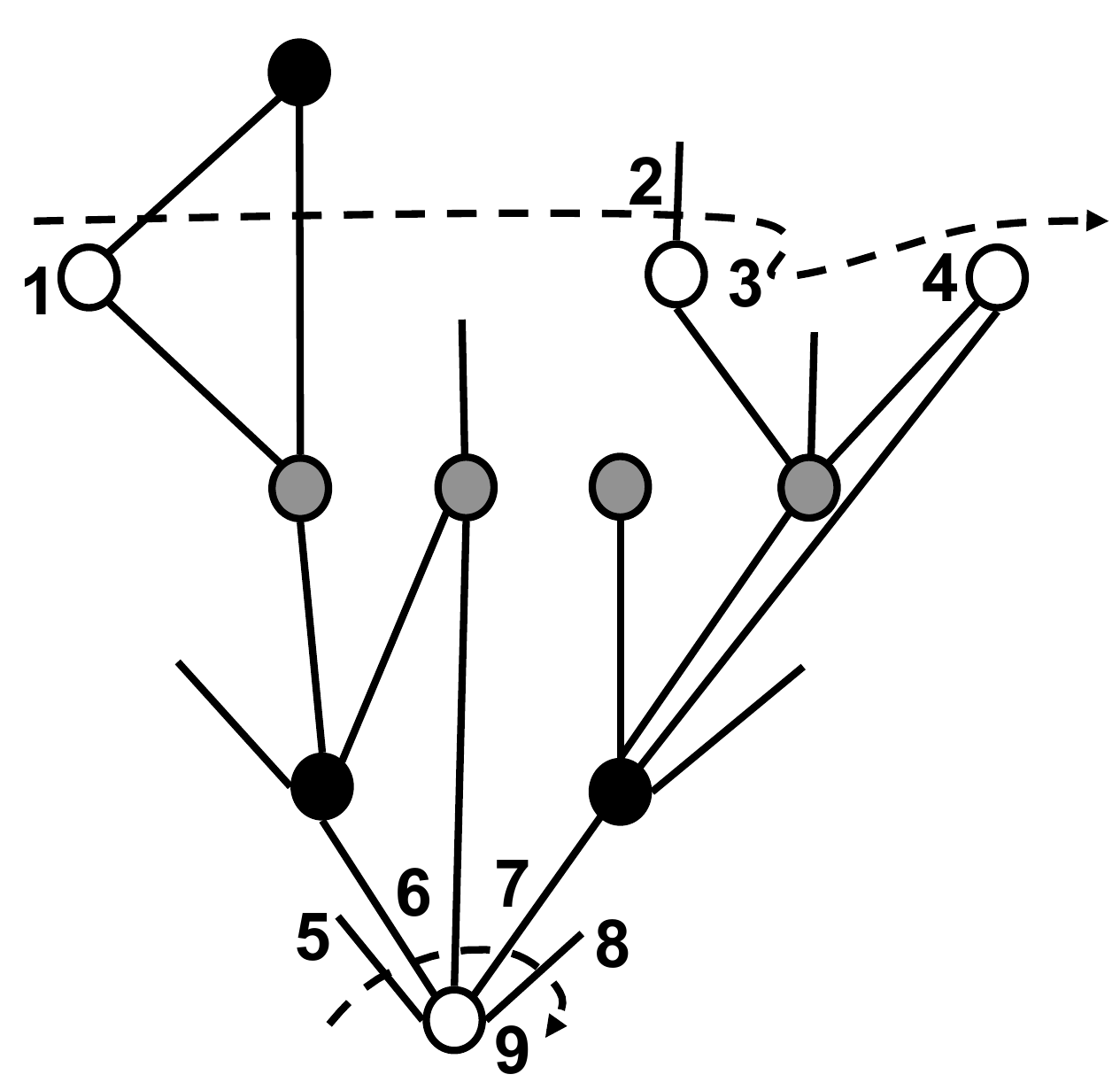}
}
\subfigure[]{
\includegraphics[width=50mm]{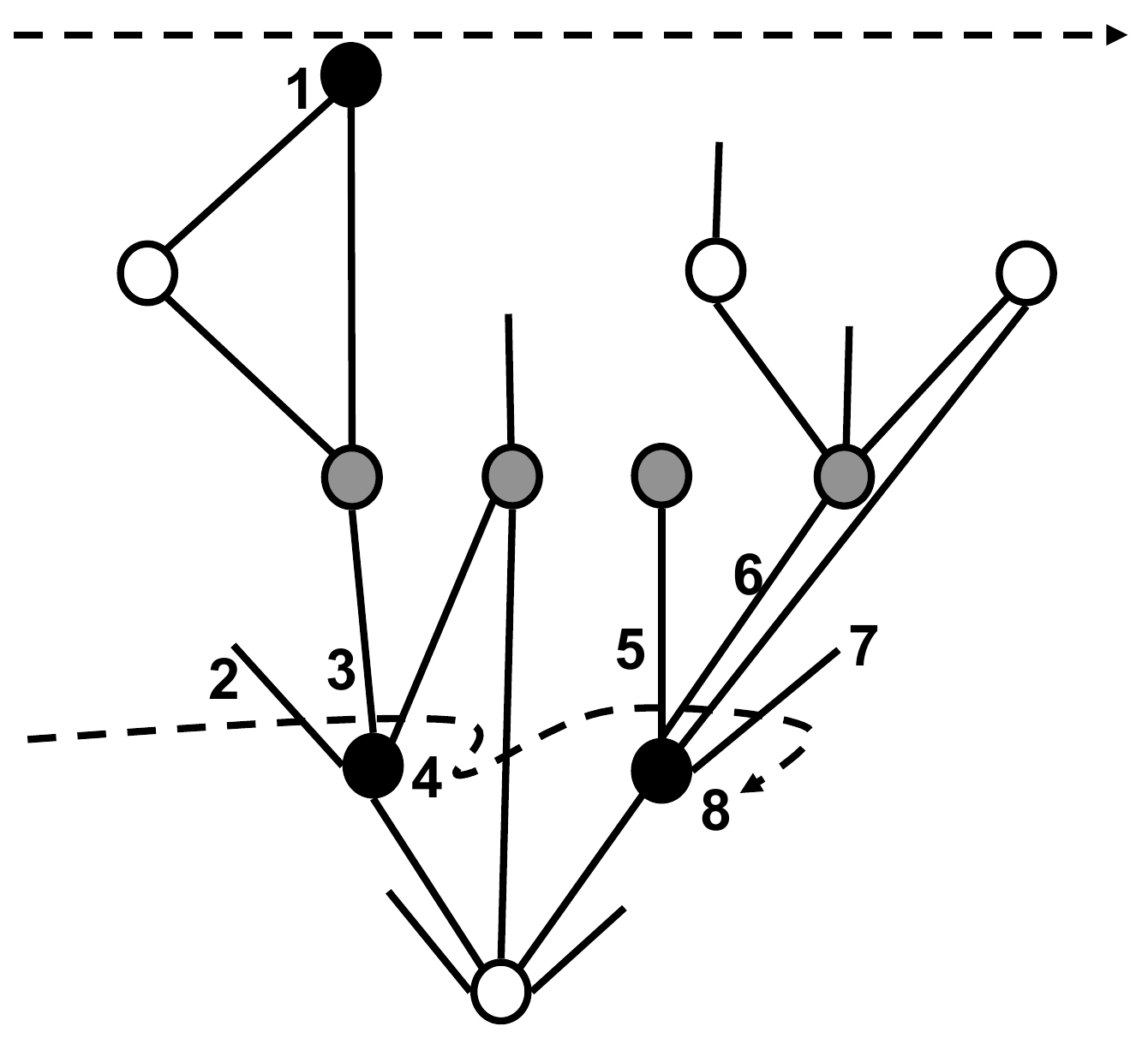}
}
\subfigure[]{
\includegraphics[width=50mm]{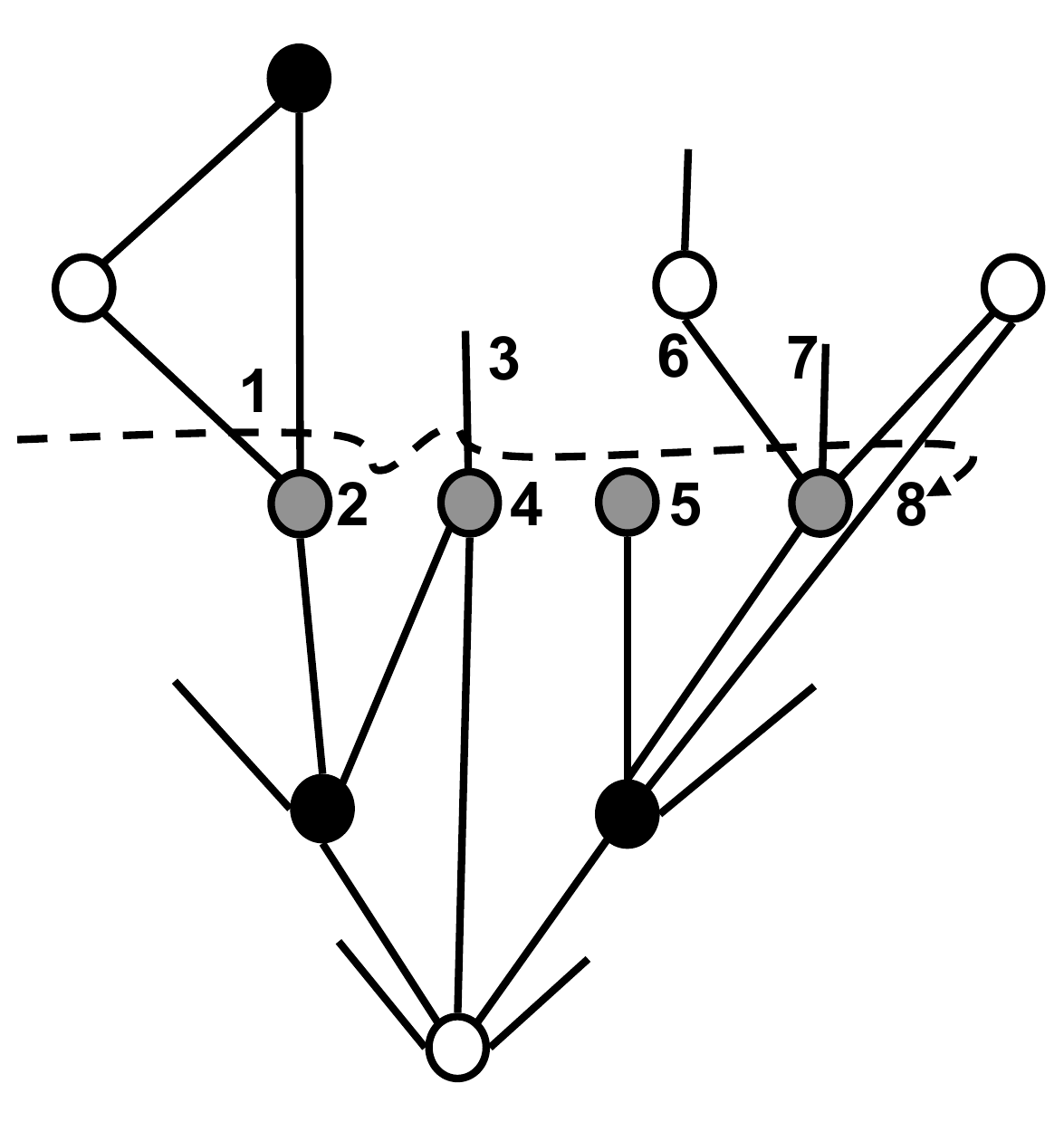}
}
\caption{Examples of (a) White, (b) black, and (c) grey reverse level traversals (RLT) on the cactus tree of Example \ref{ex:cactree}.}
\label{fig:RLT}
\end{figure}

\begin{defn}[reverse level traversals (RLT)]
For trees $\widetilde{\tau} \in \widetilde{\mathcal{CT}}(\la,\mu,\nu,g,w,b)$ (or for $\tau \in \mathcal{CT}(\la,\mu,\nu,g,w,b)$), we define the white {\em Reverse Levels Traversal (RLT)} as the following linear order in $\widetilde{\tau}$ of the white vertices and their children.
% (black vertices) and of the thorns connected to the white vertices:

We divide the white vertices of $\widetilde{\tau}$ into {\em levels} depending on their height from the root (where the height is defined as the number of edges in the shortest path with vertex sequence white-grey-black-white-\ldots to the root). So the first level consists of the root, the second level consists of the white vertices at height $3$ from the root, etc. 
The white RLT is a traversal of all the white vertices and their children (thorns, edges, triangles) from left to right, first at the level of maximum height, then the level of second maximum height,\ldots up to the root vertex. The children of each white vertex are traversed from left to right before the vertex itself \footnote{In words, the white RLT can be viewed as a {\em reverse breadth first traversal} of the white vertices and their children.}  

%\begin{itemize}
%\item[(i)] the children (either black vertices or thorns if any) of the leftmost white vertex of the top white level are traversed from left to right,
%
%\item[(ii)] then the leftmost white vertex of the top white level is traversed,
%
%\item[(iii)] if $i$ white vertices $(1,\ldots,i)$ (1 being the leftmost) of the $\ell$th white level ($\ell>1$) have been traversed and there is a white vertex $i+1$ at the $\ell$th level on the right of $i$, the children of $i+1$ are traversed, then vertex $i+1$ is traversed. Otherwise, the children of the leftmost white vertex of the white $(\ell-1)$th level are traversed, followed by the white vertex itself.
%
%\item[(iv)] The root vertex is the last to be traversed.
%\end{itemize}
\end{defn}
The black and grey RLT are defined similarly with respect to black vertices and their children, and grey vertices and their children. Figure \ref{fig:RLT} depicts the three RLT for the cactus tree in Example \ref{ex:cactree}.

\begin{defn}[Partial permutations]
Given two sets $X$ and $Y$ and a non negative integer $m$, let $\mathcal{PP}(X,Y,m)$ be the set of bijections from any $m$-subset of $X$ to any $m$-subset of $Y$ These bijections are called {\em partial permutations}. Then $|\mathcal{PP}(X,Y,m)| = \binom{|X|}{m} \binom{|Y|}{m} m!$. 
\end{defn}

\subsection{Bijective mapping $\Theta$ for 3-cacti that preserves type}

We proceed with the description of $\Theta$:
$$
\Theta^n_{\la,\mu,\nu}:\mathcal{C}(\la,\mu,\nu) \sr{\sim}{\to}  \widetilde{\mathcal{CT}}(\la,\mu,\nu) \times \mathfrak{S}_{n+1-\ell(\la)-\ell(\nu)+b} \times \mathfrak{S}_{n-\ell(\mu)-\ell(\nu)+w} \times \mathcal{OP}_{\ell(\nu)-w-b}^{(n+1-\ell(\la)-\ell(\mu)+g)}.
$$
 Within the construction, we use 
$$
(p_1, p_2, p_3) := (\ell(\la), \ell(\mu), \ell(\nu)).
$$

\subsubsection{The cactus tree $\widetilde{\tau}$}
Let $(\pi_1,\pi_2,\pi_3,\al_1,\al_2) \in \mathcal{C}(\la,\mu,\nu)$. We construct a cactus tree $\widetilde{\tau} \in \widetilde{\mathcal{CT}}(\la,\mu,\nu,g,w,b)$ following the procedure below.

\noindent (i) \emph{Cactus tree}: The first step is to construct the cactus tree $\tau$ and relabeling permutations in the same way as in \cite{V}. For completeness, we also include here the construction. 
%Let $\pi_i^{(j)}$ for $i=1,2,3$ and $j=1,\ldots,p_i$ be the blocks of the partitions $\pi_i$. 
Let $m_2'^{(j)}$ ($1\leq j \leq p_2$) be the maximal element of $\al_3^{-1}\al_2^{-1}(\pi_2^{(j)})$ and $m_i^{(j)}$ for $i=1,3$ ($1\leq j \leq p_i$) be the maximal element of the block $\pi_i^{(j)}$. %We assign the index $p_1$ to the block of $\pi_1$ containing $1$, and suppose that the indexing of the other blocks is arbitrary.

We first construct the labeled tricolored tree $T$ with $p_1$ white, $p_2$ black, and $p_3$ grey vertices satisfying:  the root of $T$ is the white vertex with label $p_1$ and the incidence relations and order of children are given in Table \ref{incidrules}.
\begin{table}[h]
$$
\begin{array}{|l|l|} \hline
\text{Incidence relations} & \text{Order of children}\\ \hline
\begin{array}{l}
\scriptstyle{\text{ for } 1\leq j\leq p_2} \\
\begin{tikzpicture}[thick,scale=0.8]
\node (r) at ( 0,0) [wv,label=left:$i$] {}; 
\node (bc2) at (-0.5,1) [bv,label=left:$j$] {};
\draw [-] (r) -- (bc2); 
\end{tikzpicture}\\
\text{ if } {m'}_2^{(j)} \in \al_3^{-1}\al_2^{-1}(\pi^{(i)}_1) 
\end{array}
&
\begin{array}{l}
\\
\begin{tikzpicture}[thick,scale=0.8]
\node (r) at ( 0,0) [wv] {}; 
\node (bc1) at (-0.5,1) [bv,label=left:$j_1$] {};
\node (bc2) at (0.5,1) [bv,label=right:$j_2$] {};
\draw [-] (r) -- (bc2);
\draw [-] (r) -- (bc1);
\end{tikzpicture}
\\
\text{ if } 
\al_2\al_3({m'}_2^{(j_1)}) <\al_2\al_3({m'}_2^{(j_2)})
\end{array}
\\[2mm] \hline 
\begin{array}{l}
\scriptstyle{\text{ for } 1\leq k \leq p_3}\\
\begin{tikzpicture}[thick,scale=0.8]
\node (r) at ( 0,0) [bv,label=left:$j$] {}; 
\node (gc) at (-0.5,1) [gv,label=left:$k$] {};
\draw [-] (r) -- (gc); 
\end{tikzpicture}\\
\text{ if } m_3^{(k)} \in \al_3^{-1}(\pi_2^{(j)}) 
\end{array}
&
\begin{array}{l}
\\
\begin{tikzpicture}[thick,scale=0.8]
\node (r) at ( 0,0) [bv] {}; 
\node (bc1) at (-0.5,1) [gv,label=left:$k_1$] {};
\node (bc2) at (0.5,1) [gv,label=right:$k_2$] {};
\draw [-] (r) -- (bc2);
\draw [-] (r) -- (bc1); 
\end{tikzpicture}
\\ 
\text{ if } 
\al_3^{-1}\al_2^{-1}\al_3(m_3^{(k_1)}) < \al_3^{-1}\al_2^{-1}\al_3(m_3^{(k_2)})
\end{array}
\\[2mm]  \hline 
\begin{array}{l}
\scriptstyle{\text{ for } 1\leq i \leq p_1-1} \\
\begin{tikzpicture}[thick,scale=0.8]
\node (r) at ( 0,0) [gv,label=left:$k$] {}; 
\node (gc) at (-0.5,1) [wv,label=left:$i$] {};
\draw [-] (r) -- (gc); 
\end{tikzpicture}\\
\text{ if } m_1^{(i)} \in \pi_3^{(k)} 
\end{array}
&
\begin{array}{l}
\\
\begin{tikzpicture}[thick,scale=0.8]
\node (r) at ( 0,0) [gv] {}; 
\node (bc1) at (-0.5,1) [wv,label=left:$i_1$] {};
\node (bc2) at (0.5,1) [wv,label=right:$i_2$] {};
\draw [-] (r) -- (bc2);
\draw [-] (r) -- (bc1); 
\end{tikzpicture}
 \\ 
\text{ if } 
\al_3^{-1}(m_1^{(i_1)}) < \al_3^{-1}(m_1^{(i_2)})\\
\end{array}
\\ \hline
\end{array}
$$
\caption{Incidence relations and order of children of the labeled tricolored tree $T$. Each row of the table shows the incidence relations and order of children of the white, black and grey vertices respectively. }
\label{incidrules}
\end{table}

\begin{lem}[\cite{V}]
The procedure above defines a labeled $3$-colored tree $T$.
\end{lem}

%\begin{proof}
%Assume that a non root white vertex $i_1$ is a child of grey vertex $k$, that is in its turn a child of black vertex $j$. Suppose now that $j$ is a child of white vertex $i_2$. Following our construction rules, we have
%$m_{1}^{(i_1)} \in \pi_{3}^{(k)}$, so that $m_{1}^{(i_1)} \leq m_{3}^{(k)}$. Then  $m_{3}^{(k)} \leq m_{2}^{'(j)}$
%as $m_{3}^{(k)} \in \alpha_3^{-1}(\pi_{2}^{(j)}) = \alpha_3^{-1}\circ \alpha_2^{-1}(\pi_{2}^{(j)})$($\pi_2$ is stable by $\alpha_2$). Finally, as $\alpha_2 \circ \alpha_3 (m_{2}^{'(j)}) \in \pi_1^{(i_2)}$,  $\gamma(m_{2}^{'(j)}) = \alpha_1 \circ \alpha_2 \circ \alpha_3 (m_{2}^{'(j)})\in \alpha_1(\pi_1^{(i_2)}) = \pi_1^{(i_2)}$ and $\gamma(m_{2}^{'(j)}) \leq m_{1}^{(i_2)}$. Two cases can occur:
%\begin{itemize}
%\item $m_{2}^{'(j)} \neq N$ then $\gamma(m_{2}^{'(j)}) = m_{2}^{'(j)}+1 \leq m_{1}^{(i_2)}$. As a direct consequence, $m_{1}^{(i_1)} <
%m_{1}^{(i_2)}$
%\item {\tt} $m_{2}^{'(j)} = N$ and $\gamma(m_{2}^{'(j)}) = 1 \in \pi_1^{(i_2)}$. As a direct consequence, $i_2$ is the root vertex.
%\end{itemize}
%Using a similar argument, one can show that if black vertex $j_1$ (resp. grey vertex $k_1$) is on a descending branch of black vertex $j_2$ (resp. grey vertex $k_2$), then $m_{2}^{'(j_1)} < m_{2}^{'(j_2)}$ (resp. $m_{3}^{(k_1)} < m_{3}^{(k_2)}$).
%Consequently, the white, black, and grey vertices' labels on each branch of the graph are strictly increasing and the last vertex of the branch is always the root vertex (by convention, the block containing $1$). 
%\end{proof}

We construct the labeled cactus tree $\Upsilon$ from $T$ by forming triangles children of the different vertices following the rules in Table \ref{rulestriang}.
\begin{table}[h]
$$
\begin{array} {|l|l|l|} \hline
\multicolumn{3}{|l|}{\text{Rules for adding triangles}} \\ \hline
\begin{array}{c}
\scriptstyle{\text{for } 1\leq i \leq p_1-1,\,\, 1\leq j \leq p_2}\\ 
\begin{tikzpicture}[thick,scale=0.8]
\node (r) at ( 0,0) [gv] {}; 
\node (gc) at (0,1) [wv,label=left:$i$] {};
\node (gc2) at (0.5,2) [bv,label=left:$j$] {};
\draw [-] (r) -- (gc); 
\draw [-] (r) -- (gc2);
\draw [-] (gc) -- (gc2);  
\end{tikzpicture}
\\
\text{ if  } \al_2\al_3({m'}_2^{(j)}) = m_1^{(i)}\\
\quad
\end{array}
&
\begin{array}{c}
\scriptstyle{\text{for } 1\leq j \leq p_2,\,\, 1\leq k \leq p_3}\\ 
\begin{tikzpicture}[thick,scale=0.8]
\node (r) at ( 0,0) [wv] {}; 
\node (gc) at (0,1) [bv,label=left:$j$] {};
\node (gc2) at (0.5,2) [gv,label=left:$k$] {};
\draw [-] (r) -- (gc); 
\draw [-] (r) -- (gc2);
\draw [-] (gc) -- (gc2); 
\end{tikzpicture}
\\
\text{ if } \al_3^{-1}\al_2^{-1}\al_3(m_3^{(k)}) = {m'}_2^{(j)} \\
(\text{i.e.}\,\, \al_3(m_3^{(k)}) = \al_2\al_3({m'}^{(j)}_2) )
\end{array}
&
\begin{array}{c}
\scriptstyle{\text{for } 1\leq k \leq p_3,\,\, 1\leq i \leq p_1-1}\\ 
\begin{tikzpicture}[thick,scale=0.8]
\node (r) at ( 0,0) [bv] {}; 
\node (gc) at (0,1) [gv,label=left:$k$] {};
\node (gc2) at (0.5,2) [wv,label=left:$i$] {};
\draw [-] (r) -- (gc); 
\draw [-] (r) -- (gc2);
\draw [-] (gc) -- (gc2); 
\end{tikzpicture}
\\
\text{ if } \al_3^{-1}(m_1^{(i)})=m_3^{(k)} \\
(\text{i.e.}\,\, m_1^{(i)} = \al_3(m_3^{(k)}) )
\end{array} 
\\ \hline
\end{array}
$$
\caption{Rules for forming triangles children of white, black and grey vertices in the labeled tricolored tree $T$ in order to obtain the labeled tricolored cactus tree $\Upsilon$.}
\label{rulestriang}
\end{table}

Finally, we remove the labels of $\Upsilon$ to obtain a cactus tree $\tau$. 

\begin{figure}[h]
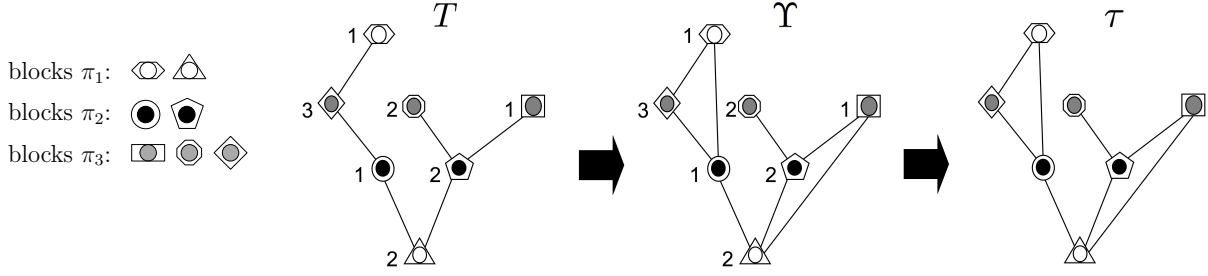

\begin{minipage}[b]{0.2\linewidth}
\includegraphics[height=3.2cm]{labelsvertices}
\end{minipage}
\begin{minipage}[b]{0.8\linewidth}
\includegraphics[clip = true, width=130mm, trim=0mm 0mm 0mm 0mm]{TGAmmaTau}
\end{minipage}
\caption{tricolored tree $T$ and cactus trees $\Upsilon$ and $\tau$ associated to Example  \ref{ex: pc1}. The tree $T$ is the $3$-colored labeled tree built following rules in Table~\ref{incidrules}. The labeled cactus tree $\Upsilon$ is obtained from $T$ by forming triangles following the rules in Table~\ref{rulestriang}. The cactus tree $\tau$ is obtained from $\Upsilon'$ by removing the labels of the vertices.}
\label{bij : TGammaTau} 
\end{figure}

\begin{exm}
\label{bij : ex}
The construction of $T$, $\Upsilon$, and $\tau$ for the partitioned cactus in Example \ref{ex: pc1} is depicted in Figure \ref{bij : TGammaTau}.
\end{exm}

\noindent (ii) \emph{Relabeling permutations}: These permutations $\theta_1,\theta_2$, and $\theta_3$ are defined by considering the {\em reverse}  labeled cactus tree $\Upsilon'$ resulting from the labeling of $\tau$, based on three independent reverse-labeling procedure for white, black, and grey vertices. We do a white RLT and label the {\em white vertices only} (as they are traversed) with the labels $1,2,\ldots, p_1$. Similarly, we do a black RLT and label the {\em black vertices only} with labels $1,2,\ldots, p_2$, and do a grey RLT to label the {\em grey vertices only} with labels $1,2,\ldots, p_3$. Next, we reindex the blocks of $\pi_1, \pi_2$ and $\pi_3$ using the new indices from $\Upsilon'$: if a white vertex is labeled $i$ in $T$ and $i'$ in $\Upsilon'$, we set $\pi_1^{i'}=\pi_1^{(i)}$ (and $m_1^{i'}=m_1^{(i)}$). Black and grey blocks are reindexed in a similar fashion. Let $u^i,v^j,w^k$ be the strings obtained by writing the elements of $\pi_1^i,\pi_2^j,\pi_3^k$ respectively in increasing order. Denote $u=u^1\cdots u^{p_1}$, $v=v^1\cdots v^{p_2}$, $w=w^1\cdots w^{p_3}$ the concatenations of the strings defined above. We define $\theta_1 \in \mathfrak{S}_n$ by setting $u$ as the first line and $[n]$ as the second line of the two-line representation of this permutation. Similarly, we define the relabeling permutations $\theta_2$ and $\theta_3$ from $v$ and $w$ respectively. 

\begin{figure}[h]
\begin{minipage}[b]{0.18\linewidth}
\includegraphics[height=2.8cm]{labelsvertices}
\end{minipage}
\begin{minipage}[b]{0.8\linewidth}
\includegraphics[clip = true, width=130mm, trim=0mm 0mm 0mm 0mm]{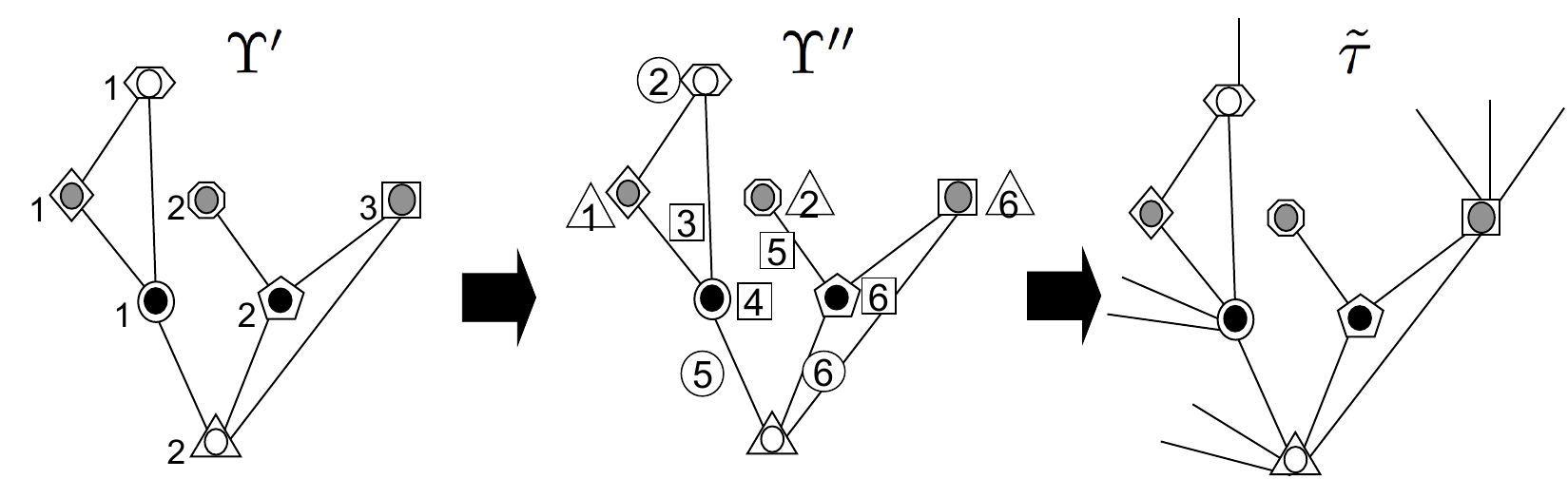}
\end{minipage}
\caption{Construction of the cactus tree $\widetilde{\tau}$ associated to Example \ref{ex: pc1}}
\label{bij : GammaPGammaSTauTilde} 
\end{figure}

\begin{exm}
\label{pperm : ex} 
Following up on Example \ref{bij : ex}, we construct the relabeling permutations $\theta_1,\theta_2$, and $\theta_3$. We have:

\noindent $\theta_{1}  = \left(\begin{array}{cc|cccc}
4 & 5 & 1 & 2 & 3 & 6\\
1 & 2 & 3 & 4 & 5 & 6
\end{array}\right),\,\,
\theta_{2}  = \left(\begin{array}{cccc|cc}
1 & 3 & 4 & 5 & 2 & 6\\
1 & 2 & 3 & 4 & 5 & 6
\end{array}\right), \,\,
\theta_3 = \left(
\begin{array}{c|c|cccc}
5 & 2 & 1 & 3 & 4 & 6\\
1 & 2 & 3 & 4 & 5 & 6
\end{array}
\right).
$

\end{exm}
 We define the multisets:
\begin{eqnarray*}
S_1&=&{\Big \{}\theta_1(m_1^i){\Big \}}_{i=1}^{p_1-1} \cup {\Big \{}\theta_1(\al_2\al_3({m'}_2^j)){\Big\}}_{j=1}^{p_2} \subset [n],\\
S_2&=& {\Big \{}\theta_2({m'}_2^j){\Big\}}_{j=1}^{p_2-1} \cup {\Big\{}\theta_2(\al_3^{-1}\al_2^{-1}\al_3(m_3^k)){\Big\}}_{k=1}^{p_3} \subset [n-1],\\ 
S_3&=&{\Big\{}\theta_3(m_3^k){\Big\}}_{k=1}^{p_3-1} \cup {\Big\{}\theta_3(\al_3^{-1}(m_1^i)){\Big\}}_{i=1}^{p_1-1} \subset [n-1].
\end{eqnarray*}
They are multisets since we allow some elements to be repeated once when there are triangles in $\Upsilon'$.
Note that the sizes of the underlying sets of $S_1,S_2$, and $S_3$ are $p_1+p_2-1-g,\, p_2+p_3-1-w$, and $p_1+p_3-2-b$ respectively. We use these multisets to label the vertices, the edges and triangles of the cactus tree $\tau$ with three types of labels: {\em circle} for $\theta_1$, {\em square} for $\theta_2$, and {\em triangle} for $\theta_3$ (represented as labels \includegraphics[height=0.35cm]{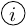}, \includegraphics[height=0.35cm]{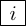}, and : \includegraphics[height=0.35cm]{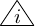} as illustrated in Figure~\ref{bij : GammaPGammaSTauTilde}). We assign $\theta_1(m_1^i)$ to the white vertices indexed $i$ in $\Upsilon'$, $\theta_2({m'}^j_2)$ to the black vertices indexed $j$ in $\Upsilon'$, and $\theta_3(m_3^k)$ to the grey vertices indexed $k$ in $\Upsilon'$. Children of a white vertex (edges and triangles) are labeled $\theta_1(\al_2\al_3({m'}_2^j))$ if the black vertex at the end of the edge or within the triangle is indexed by $j$ in $\Upsilon'$. The children of the black and grey vertices are labeled in a similar fashion with ${\big \{}\theta_2(\al_3^{-1}\al_2^{-1}\al_3(m_3^k)){\big \}}_{k=1}^{p_3}$ and ${\big \{}\theta_3(\al_3^{-1}(m_1^i)){\big \}}_{i=1}^{p_1-1}$ respectively.

Let $\Upsilon''$ be the resulting cactus tree with these new additional labels. 

%Note that labels (of the same type) of a vertex and its child are equal whenever the two vertices belong to a triangle rooted in another vertex. 
 %picture

Let $\ov{S_i}$ (for $i=1,2,3$) be the ordered multiset obtained by arranging the elements of $S_i$ in non-decreasing order (for $i=1,2,3$).

\begin{lem}
\label{d}
Let $\mathbf{d}=(d_1,\ldots,d_{p_1+p_2-1})$ ($\mathbf{d}'=(d'_1,\ldots,d'_{p_2+p_3})$, and $\mathbf{d}''=(d''_1,\ldots,d''_{p_1+p_3-1})$ resp.) be the ordered set of labels of the first type (second and third type resp.) obtained by traversing $\Upsilon''$ up to, but not including the root vertex according to the white RLT (black and grey RLT resp.) defined in Section 3. We have:
$$
\mathbf{d}=\ov{S_1},\quad \mathbf{d}'=\ov{S_2}, \quad \mathbf{d}''=\ov{S_3}.
$$
\end{lem}
 
\begin{proof}
Let $\theta_1(m_1^0)=0$. By construction, if $\theta_1(\al_2\al_3({m'}_2^j))$ in $\Upsilon''$ is the label of a child of a white vertex with label $\theta_1(m_1^i)$, $1\leq i \leq p_1$, then $\al_2\al_3({m'}_2^j) \in \pi_1^i$ and $\al_2\al_3({m'}_2^j)\leq m_1^i$. As $\theta_1$ is increasing among the blocks and within them then: 
$$
\theta_1(m_1^{i-1}) < \theta_1(\al_2\al_3({m'}_2^j))<\theta_1(m_1^i),
$$ 
%where the latter is an equality when the white vertex $i$ is the middle vertex of a triangle with black vertex $j$ and rooted at the grey ancestor of vertex $i$. 
Now, if two children of the same white vertex have respective labels $\theta_1(\al_2\al_3({m'}_2^{j_1}))$ and $\theta_1(\al_2\al_3({m'}_2^{j_2}))$ in $\Upsilon''$, and black vertex $j_1$ is to the left of $j_2$; then by construction we have $\al_2\al_3({m'}_2^{j_1})<\al_2\al_3({m'}_2^{j_2})$. Again, since $\theta_1$ is increasing within blocks of $\pi_1$ then $\theta_1(\al_2\al_3({m'}_2^{j_1}))<\theta_1(\al_2\al_3({m'}_2^{j_2}))$. Finally, the white RLT of the circle labels in $\Upsilon''$ (up to but not including the root) yields $\ov{S_1}$. Similarly, black and grey RLTs of the square and triangle labels in $\Upsilon''$ yield $\ov{S_2}$ and $\ov{S_3}$ respectively.
\end{proof}

\medskip

\noindent (iii) \emph{Thorns}: Recall the definition of $\mathbf{d},\mathbf{d}',$ and $\mathbf{d}''$ in Lemma \ref{d}. We add thorns to the white vertices for each missing element of $[n]$ in $\mathbf{d}$, and add thorns to black and grey vertices for  each missing element of $[n-1]$ in $\mathbf{d}'$ and $\mathbf{d}''$, respectively. More specifically, we add $n+1-p_1-p_2+g$ thorns to the white vertices, $n-p_2-p_3+w$ thorns to the black vertices, and $n+1-p_1-p_3+b$ thorns to the grey vertices of $\Upsilon''$ in the following fashion:
 %We know that the vectors $\mathbf{d},\mathbf{d}'$, and $\mathbf{d}''$ are increasing subsequences obtained by the RLTs of $\Upsilon''$. %Below, we describe the rules for adding thorns to the white vertices. By convention, if the vertex is a root of a triangle we do not add thorns {\em inside} it. 

\begin{itemize}
\item[1.] If $d_1 > 1$ and $d_1$ is the label of a (white) vertex, we connect $d_1-1$ thorns to it. If a child of a white vertex has label $d_1$, we connect $d_1-1$ thorns to the ascending white vertex on the left of child $d_1$.

\item[2.] For $1<l<p_1+p_2-1$, if $d_l > d_{l-1}+1$ we follow one of the four following cases:\\
(a) $d_l$ and $d_{l-1}$ are both the label of white vertices in $\Upsilon''$, white vertex $d_l$ (short for vertex corresponding to $d_l$) has no  child and it is the white vertex following $d_{l-1}$ in the white RLT of $\Upsilon''$. If so, we connect $d_l - d_{l-1}-1$ thorns to $d_{l}$.\\
(b) $d_l$ is the first label of a child and $d_{l-1}$ is the first label of a white one, then $d_{l}$ is the leftmost child of the white vertex following  $d_{l-1}$. If so, we connect $d_l - d_{l-1}-1$ thorns to the ascending white vertex of $d_l$ on its left\\
(c) $d_l$ is the first label of a white vertex and $d_{l-1}$ is the first label of a child, then $d_{l-1}$ is the rightmost child of $d_l$. If so, we connect $d_l - d_{l-1}-1$ thorns to $d_l$ on the right of $d_{l-1}$\\
(d) Finally, if $d_l$ and $d_{l-1}$ are both the first label of children, they have the same white ascending vertex. We connect $d_l - d_{l-1}-1$ thorns to the ascending white vertex between them.
\item[3.] If $d_{p_1+p_2-1}<n$, we connect $n- d_{p_1+p_2-1}-1$ thorns to the root vertex on the right of its rightmost child.
\end{itemize}

Again, we can think of this as adding a thorn to the the white vertices for each element of $[n]$ not included in $\mathbf{d}$.

A similar construction is applied to add thorns to the black and grey vertices following the sequence of integers $\mathbf{d}'$ and $\mathbf{d}''$. Finally we remove all the labels to get the cactus tree $\widetilde{\tau}$.

\begin{exm}
Figure \ref{bij : GammaPGammaSTauTilde} depicts the construction of the cactus tree $\widetilde{\tau}$ corresponding to the partitioned cactus in Example \ref{ex: pc1}.
\end{exm}

The next two lemmas show that $\widetilde{\tau}$ preserves the type of the partitioned cacti, and that $\Upsilon''$ can be recovered from $\widetilde{\tau}$ via white, black, and grey RLTs. 

\begin{lem}
\label{tau}
$\widetilde{\tau}$ as defined above belongs to $\widetilde{\mathcal{CT}}(\la,\mu,\nu,g,w,b)$ where $g,w,b$ are the number of triangles in $\widetilde{\tau}$ rooted in grey, white, and black vertices respectively.
\end{lem}

\begin{proof}
%First we check that the thorns in $\widetilde{\tau}$ are not inside the triangles. In the case of the thorns added to the a white vertex $i$, we only need to consider the instance when $i$ is the middle vertex of a triangle rooted on a grey vertex. This could apply to cases 2. (b) and (c) above. In the former, the thorns cannot be inside a triangle since the non-tree edge would be on the right of the black vertex labeled $d_l$ (recall that in a triangle, the vertices are ordered clockwise as white-black-grey-white). In the latter case, since the label $d_l <d_{l-1}$, then $d_{l}$ is not the middle vertex of a triangle. A similar conclusion is reached by looking at the rules for adding thorns to black and grey vertices.

We check the vertex degrees of $\widetilde{\tau}$. If we take two successive white vertices $i-1$ and $i$ according to  white RLT of $\Upsilon''$ with labels $\theta_1(m_1^{i-1})$ and $\theta_1(m_1^{i}), (i<p_1)$, a thorn is connected to $i$ for any missing integer of the interval $[\theta_1(m_1^{i-1}),\theta_1(m_1^{i})-1]$ in $\mathbf{d}$. This number of missing integers is equal to $\theta_1(m_1^{i})-1-\theta_1(m_1^{i-1})-f_i$ where $f_i$ is the number of children of $i$. As $i$ is not the root vertex, there is an edge between $i$ and its ancestor so that the resulting degree $deg$ for $i$ (as defined in \eqref{degdef}) is:
\begin{equation}
deg(i) = f_i+\left(\theta_1(m_1^{i})-1-\theta_1(m_1^{i-1}) - f_i \right) +1 =\theta_1(m_1^{i})-\theta_1(m_1^{i-1}), \,\,\,\, \forall i \in [p_1-1],
\end{equation}
Furthermore, $n-\theta_1(m_1^{p_1-1})-f_p$ thorns are connected to the root vertex (since $\varepsilon_{p_1}=0$) so that:
\begin{equation}
deg(p_1) = n-\theta_1(m_1^{p_1-1})
\end{equation}
But, according to the construction of $\theta$, 
\begin{eqnarray}
\theta_1(\pi_1^1) &=& [\theta_1(m_1^{1})]\\
\theta_1(\pi_1^i) &=& [\theta_1(m_1^{i})]\setminus[\theta_1(m_1^{i-1})], \; (2\leq i \leq p_1-1)\\
\theta_1(\pi_1^{p_1}) &=& [n]\setminus[\theta_1(m_1^{p_1-1})]
\end{eqnarray}
Subsequently:
\begin{equation}
deg(i) = |\pi_1^i|, \,\, \forall i \in [p_1].
\end{equation}
And $\la = type(\pi)$ is the white vertex degree distribution of $\widetilde{\tau}$. In a similar fashion, $\mu$ and $\nu$ are the black and grey vertex degree distribution of $\widetilde{\tau}$.
\end{proof}

\begin{lem}
\label{labels}
Assign {\em circle} labels $1,2,\ldots,n$ to the white vertices and their children (including thorns) in $\widetilde{\tau}$ in increasing order according to the white RLT, add two other sets of labels $1,2,\ldots,n$ ({\em square} and {\em triangle}) to the black and grey vertices and their children  in increasing order according to the black and grey RLT. The labeling of the vertices and children that are not thorns is the same as in $\Upsilon''$. 
\end{lem}

\begin{proof}
According to the construction of $\widetilde{\tau}$, we add thorns to $\Upsilon''$ when integers are missing in its RLTs so that the thorns would take these missing integers as labels when traversing the cactus tree. As a result, the labels of the vertices in the  RLTs of $\widetilde{\tau}$ are still $\mathbf{d},\mathbf{d}'$, and $\mathbf{d}''$ and since they still appear in the same order, we have the desired result. 
\end{proof}

\subsubsection{The permutations $\sigma_1$ and $\sigma_2$ and the ordered set $\chi$ and the }

In the previous subsection we explained how to obtain the cactus tree $\widetilde{\tau}$ from the partitioned $3$-cactus in $\mathcal{C}(\lambda,\mu,\nu)$. We now move on to explain how to obtain the permutations $\sigma_1$ in $\mathfrak{S}_{n+1-\ell(\lambda)-\ell(\nu)+b}$ and $\sigma_2$ in $\mathfrak{S}_{n-\ell(\lambda)-\ell(\nu)+b}$, and the ordered set $\chi$ in $\mathcal{OP}^{(n+1-\ell(\lambda)-\ell(\mu)+g)}_{\ell(\nu)-w-b}$.
\begin{itemize}
\item[(i)] {\em Permutations $\sigma_1,\sigma_2$:} Let $E$ and $F$ be the following sets: 
\begin{align*}
E &= [n]\left\backslash \, \left( {\big \{}\theta_1(m_1^i){\big \}}_{i=1}^{p_1-1} \cup {\big \{}\theta_1(\al_3(m_3^k)){\big \}}_{k=1}^{p_3-1} \right)\right.,\\
F &= [n]\left\backslash \, \left({\big \{}\theta_1(\al_2\al_3({m'}_2^j){\big \}}_{j=1}^{p_2} \cup {\big \{}\theta_1(\al_3(m_3^k)){\big \}}_{k=1}^{p_3-1} \right)\right. .
\end{align*} 
We define partial permutations $\tilde{\sigma_1}$ and $\tilde{\sigma_2}$ in the following way:
$$
\begin{array}{cccc}
\tilde{\sigma_1}: &   E &\to& [n-1]\backslash S_3\\
 & u &\mapsto& \theta_3\alpha_3^{-1}\theta_1^{-1}(u)\\
\tilde{\sigma_2}:&   F &\to& [n-1]\backslash S_2\\
 & u &\mapsto& \theta_2 \alpha_3^{-1}\alpha_2^{-1}\theta_1^{-1}(u).
\end{array}
$$
Let $\sigma_1 \in \mathfrak{S}_{n+1-p_1-p_3+b}$ and $\sigma_2 \in \mathfrak{S}_{n-p_2-p_3+w}$ be the order isomorphic permutations corresponding to $\tilde{\sigma_1}$ and $\tilde{\sigma_2}$ respectively. 

\item[(ii)] {\em Ordered set $\chi$:} We define the ordered set $\widetilde{\chi}={\big \{}\theta_1(\al_3(m_3^k)) \mid \theta_1(\al_3(m_3^k)) \notin S_1 {\big \}}_{k=1}^{p_3}$. Then, let $\rho : [n]\setminus S_1 \vdash [n- |S_1|]$ be the indexing permutation associating to any integer $i \in [n]\setminus S_1$ its position in $\ov{[n]\setminus S_1}$ where $\ov{[n]\setminus S_1}$ is the ordered (increasing) set of $[n]\setminus S_1$. The ordered set $\chi$ is defined as follows:
\begin{equation}
\chi = \rho(\widetilde{\chi})
\end{equation}
As $|S_1| = n-(\ell(\la)-1)-\ell(\mu)+g$ and $\left|{\big \{}\theta_1(\al_3(m_3^k)) \mid \theta_1(\al_3(m_3^k)) \notin S_1 {\big \}}\right| = \ell(\nu)-w-b$, $\chi$ belongs to the set $\mathcal{OP}_{\ell(\nu)-w-b}^{(n+1-\ell(\la)-\ell(\mu)+g)}$.
\end{itemize}

%$$
%E_0\backslash (E_0\cap S_1)=\{\theta_1(m_1^{i_1}),\theta_1(m_1^{i_2}),\ldots,\theta_1(m_1^{i_{p_3-w-b+1}}) \},
%$$
%where $1\leq i_1<i_2<\ldots<i_{p_3-w-b+1}\leq p_3$. Let $\chi$ be the ordered set $(\theta_1(m_1^{i_1}),\theta_1(m_1^{i_2}),\ldots,\theta_1(m_1^{i_{p_3-w-b+1}}))$.

\begin{figure}
\begin{minipage}[b]{0.08\linewidth}
$\Theta^n_{\lambda,\mu,\nu}:$
\vspace{0.4in}
\end{minipage}
\begin{minipage}[b]{0.33\linewidth}
\includegraphics[width=50mm]{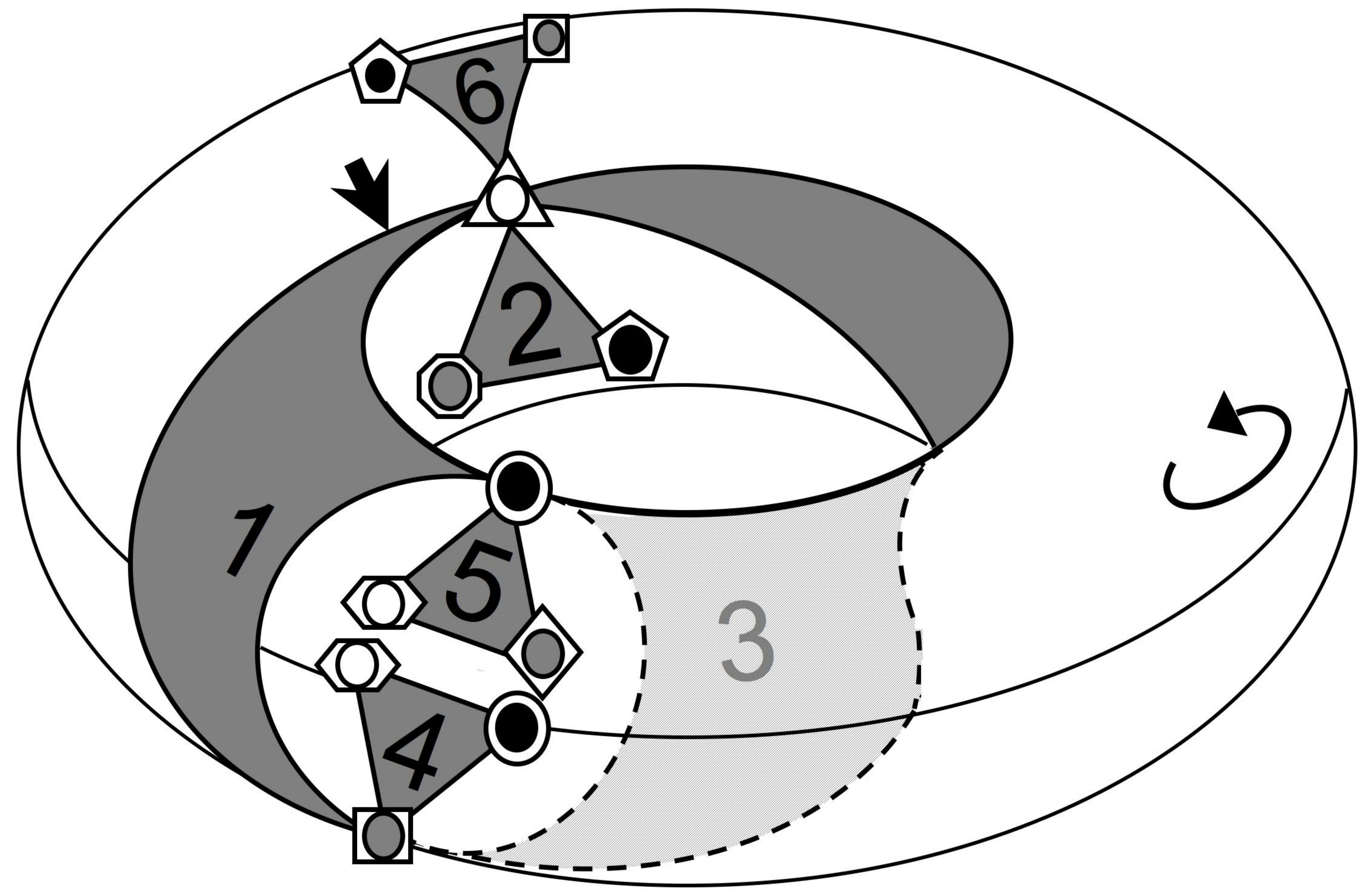}
\end{minipage}
\begin{minipage}[b]{0.56\linewidth}
\includegraphics[width=90mm]{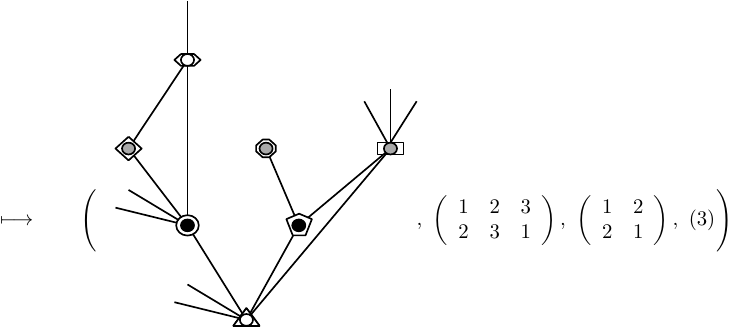}
\end{minipage}
\caption{Summary of ouput $(\widetilde{\tau},\sigma_1,\sigma_2,\chi)$ of the map $\Theta^n_{\lambda,\mu,\nu}$ applied to the partitioned $3$-cactus from Example~\ref{ex: pc1}.}
\label{fig:summarybij}
\end{figure}

\begin{exm}
Getting back to Example \ref{ex: pc1}, computing the partial permutations leads to:
$$\tilde{\sigma_1} = \left (\begin{array}{ccclccc}1&3&5\\4&5&3\end{array} \right),\,\,\,\,\,\,\tilde{\sigma_2} = \left (\begin{array}{cclcc}3&1\\1&2\end{array} \right),$$
and $$\sigma_1 = \left (\begin{array}{ccclccc}1&2&3\\2&3&1\end{array} \right),\,\,\,\,\,\,{\sigma_2} = \left (\begin{array}{cclcc}1&2\\2&1\end{array} \right).$$
For the ordered set we have: $\widetilde{\chi} = (4)$ and $\chi = (3)$.\\ 
\noindent In summary, the map $\Theta^n_{\lambda,\mu,\nu}$ applied to $(\pi_1,\pi_2,\pi_3,\alpha_1,\alpha_2)$ in $\mathcal{C}({[2^1,4^1],[2^1,4^1], [1^2,4^1]})$ from Example~\ref{ex: pc1} gives the 4-tuple $(\widetilde{\tau}, \sigma_1,\sigma_2, \chi)$ depicted in Figure~\ref{fig:summarybij}.
\end{exm}

% pending picture

%There is a bijection $\Theta^n_{\la,\mu,\nu}$ between partitioned 3-cacti $\mathcal{C}(\la,\mu,\nu)$ and tuples $(\widetilde{\tau},\sigma_1,\sigma_2,\chi)$ where $\widetilde{\tau} \in \widetilde{CT}(\la,\mu,\nu,g,w,b)$, $\sigma_1 \in S_{n+1-\ell(\la)-\ell(\nu)+b}$, $\sigma_2 \in S_{n-\ell(\mu)-\ell(\nu)+w}$, and $\chi \in \mathcal{OP}_{\ell(\nu)-w-b}(n+1-\ell(\la)-\ell(\mu)+g)$ for some $g,w,b\geq 0$.

\subsection{Showing the map $\Theta$ is a bijection}
To show that $\Theta^n_{\la,\mu,\nu}$ is a one-to-one correspondence we take any element $(\widetilde{\tau},\sigma_1,\sigma_2,\chi)$ in 
$$
\widetilde{\mathcal{CT}}(\la,\mu,\nu,g,w,b) \times \mathfrak{S}_{n+1-\ell(\la)-\ell(\nu)+b} \times \mathfrak{S}_{n-\ell(\mu)-\ell(\nu)+w} \times \mathcal{OP}_{\ell(\nu)-w-b}^{(n+1-\ell(\la)-\ell(\mu)+g)}
$$
and show that there is a unique element $(\pi_1,\pi_2, \pi_3, \alpha_1, \alpha_2)$ in $\mathcal{C}(\la,\mu,\nu)$ such that\linebreak $\Theta^n_{\la,\mu,\nu}(\pi_1,\pi_2, \pi_3, \alpha_1, \alpha_2) =(\widetilde{\tau},\sigma_1,\sigma_2,\chi)$. Let $p_1=\ell(\la)$, $p_2 = \ell(\mu)$, and $p_3=\ell(\nu)$. We proceed with a two step proof:
\begin{itemize} 

\item[(i)] The first step is to notice that $(\widetilde{\tau},\sigma_1,\sigma_2,\chi)$ defines a unique cactus tree $\tau$ belonging to $\mathcal{CT}(p_{1}, p_{2}, p_{3}, g, w, b)$, unique multisets $\{S_i\}_{1 \leq i \leq 3}$, as well as a unique ordered set $\widetilde{\chi}$ belonging to $\mathcal{OP}_{p_3-w-b}^{(n+1-p_1-p_2+g)}$. Labeling each vertex and children of $\tau$ with $1,2,\ldots,n$ in increasing order according to the three reverse levels traversals and removing the three sets of thorns (together with their labels) gives a labeled cactus tree $\Upsilon''$ that leads to $\widetilde{\tau}$ according to $\Theta$. This labeled cactus tree is the unique one that can lead to $\widetilde{\tau}$ since within $\Theta$, $\widetilde{\tau}$, and $\Upsilon''$ have the same underlying cactus tree structure $\tau$, and according to Lemma \ref{labels}, $\widetilde{\tau}$ determines the labels of $\Upsilon''$.\\
Then, using Lemma \ref{d}, the three series of labels (except the root's) in $\Upsilon''$ are by construction the three sets $\{S_i\}_{1 \leq i \leq 3}$. The knowledge of $S_1$ and $\chi$ uniquely determines $\widetilde{\chi}$. As a result, exactly one $7$-tuple $(\tau, S_1, S_2, S_3, \sigma_1, \sigma_2, \widetilde{\chi})$ is associated to $(\widetilde{\tau},\sigma_1,\sigma_2,\chi)$ by the final steps of the mapping $\Theta$.

\item[(ii)] The bijection $\Theta_{n,p_1,p_2,p_3}$ in \cite{V} is identical to the first steps (up to the construction of  $\tau, S_1, S_2, S_3, \sigma_1, \sigma_2$ and $\widetilde{\chi}$) of $\Theta^n_{\la,\mu,\nu}$. Therefore by \cite[Sec. 6]{V} there is a unique $5$-tuble $(\pi_1,\pi_2,\pi_3,\alpha_1,\alpha_2)$ in $\mathcal{C}(p_1,p_2,p_3,n)=\bigcup_{\ell(\la)=p_1,\ell(\mu)=p_2,\ell(\nu) = p_3}\mathcal{C}(\la,\mu,\nu)$ mapped to the $7$-tuple $(\tau, S_1, S_2, S_3, \sigma_1, \sigma_2, \widetilde{\chi})$ by $\Theta_{n,p_1,p_2,p_3}$ and equivalently by the first steps of $\Theta^n_{\la,\mu,\nu}$. \\
According to \cite{V}, the  types of $\pi_1$, $\pi_2$, and $\pi_3$ are directly recovered from $\{S_i\}_{1 \leq i \leq 3}$ and $\tau$. 
Furthermore, using Lemma \ref{tau}, the vertex degree distribution of $\widetilde{\tau}$ is equal to the type of the partitions encoded by the elements in $\{S_i\}_{1 \leq i \leq 3}$ corresponding to the relabeling of the maximum elements of the blocks. Finally, as the vertex degree distribution in $\widetilde{\tau}$ is $(\la,\mu,\nu)$, so is the type of $(\pi_1,\pi_2,\pi_3)$. Therefore, $(\pi_1,\pi_2,\pi_3,\alpha_1,\alpha_2)$ belongs to $\mathcal{C}(\la,\mu,\nu)$ as desired.

\end{itemize}

%\subsection{Connection to combinatorial results about spectra of normally distributed random matrices}
%

%
%\begin{rem}
%If we use $p_{\la_i}(X)=trace(X^{\la_i})$ for a square matrix $X$, then by \cite{HSS} and \cite{GJ95} we can view the RHS of Equation \eqref{eq2} as $\mathcal{E}_{U}(p_{n}(XUYU^*))$ where the expectation is over $n\times n$ matrices $U$ whose entries are independent standard normal random complex variables and $X,Y$ are arbitrary fixed Hermitian complex matrices with eigenvalues marked by the indeterminates $x_1,x_2,\ldots$ and $y_1,y_2,\ldots$ are the eigenvalues of $X$ and $Y$ respectively. Similarly, by linearity of expectation and the identity $k_{\la^{(1)}\la^{(2)}\la^{(3)}}^{\mu} = \sum_{\nu \vdash n}k_{\la^{(1)}\la^{(2)}}^{\nu}k_{\nu\la^{(3)}}^{\mu}$, the RHS of Equation \eqref{eq1} can be viewed as $\mathcal{E}_{U}(\mathcal{E}_V(p_{n}(XUYU^*VZV^*)))$ where the expectation is over pairs of $n\times n$ matrices $U$ and $V$ whose entries are independent standard normal random complex variables.  
%\end{rem}

%\subsection{Appendix: Derivation of number of cactus trees} 

\section{Proof  of Proposition \ref{numcact}: computation of the number of cactus trees} \label{sect:pfpropnumcact}
 
In this section we prove Proposition \ref{numcact} where we compute the cardinality of the set $\widetilde{\mathcal{CT}}(\lambda,\mu,\nu,g,w,b)$. To do this, we consider its generating function $F$: 
\begin{equation}
F = \sum_{\lambda, \mu, \nu \vdash n} \sum_{g,w,b \geq 0} |\widetilde{\mathcal{CT}}(\lambda,\mu,\nu,g,w,b)| x_1^{\ell(\lambda)}x_2^{\ell(\mu)}x_3^{\ell(\nu)}x_4^{g}x_5^{w}x_6^{b}\mathbf{t}^{\mathbf{n}(\lambda)}\mathbf{u}^{\mathbf{n}(\mu)}\mathbf{v}^{\mathbf{n}(\nu)}.
\end{equation}
That is, the white, black, and grey vertices are marked respectively by  indeterminates $x_1, x_2$ and $x_3$. Triangles children of a grey, white, and black vertex are marked respectively
by $x_4, x_5,$ and $x_6$. Furthermore, $t_i,u_j,$ and $v_k$ mark respectively white vertices of degree $i$, black vertices of degree $j$ and grey vertices of degree $k$. And $\bold{t} = (t_1,t_2,\ldots)$, $\bold{u} = (u_1,u_2,\ldots)$, $\bold{v} = (v_1,v_2,\ldots)$ and $\mathbf{n}(\epsilon) = (n_1(\epsilon),n_2(\epsilon),\ldots)$ for $\epsilon \vdash n$ where $n_i(\epsilon)$ is the number of $i$ parts of $\epsilon$. 

%As a direct consequence:
%\begin{equation}
%|\widetilde{\mathcal{CT}}(\lambda,\mu,\nu,g,w,b)| = [x_1^{\ell(\lambda)}x_2^{\ell(\mu)}x_3^{\ell(\nu)}x_4^{g}x_5^{w}x_6^{b}\mathbf{t}^{\mathbf{n}(\lambda)}\mathbf{u}^{\mathbf{n}(\mu)}\mathbf{v}^{\mathbf{n}(\nu)}]\,\,F
%\end{equation}

The evaluation of $F$ is performed thanks to the multivariate Lagrange inversion theorem (see e.g. \cite[1.2.13]{GJCE}). We propose a recursive decomposition of the desired set of cactus trees sharing similar ideas with \cite{GJ92}.\\

In a similar fashion as in \cite{GJ92}, we introduce $W$, $B$, and $G$ as the generating functions of the sets $\mathcal{W}$, $\mathcal{B}$, and $\mathcal{G}$ of non empty {\it planted} cactus trees with respectively white, black, and grey root vertices.
Construction rules for these sets of cactus trees are identical to those of $\widetilde{\mathcal{CT}}(\lambda,\mu,\nu,g,w,b)$ with the only exception that an additional {\it planted} edge is connected to the root vertex on the left of the leftmost child (vertex or thorn). We take this additional edge into account in the root's degree. Finally, let $T_g$, $T_w$, and $T_b$ be respectively the generating functions of triangles children of a grey, white, and black vertices. Immediately:
\begin{eqnarray}
T_g = x_4\\
T_w = x_5\\
T_b = x_6
\end{eqnarray} 

Any cactus tree in $\widetilde{\mathcal{CT}}(\lambda,\mu,\nu,g,w,b)$ can be decomposed in a tuple of planted cactus trees in $\mathcal{W}$, $\mathcal{B}$, and $\mathcal{G}$. The rule for the decomposition is based on the nature of the leftmost child of the white root in a given cactus tree $\tau$ of $\widetilde{\mathcal{CT}}(\lambda,\mu,\nu,g,w,b)$:
\begin{compactenum}
\item[(i)] If the leftmost child is a thorn then $\tau$ is equivalent to the cactus tree in $\mathcal{W}$ with the planted edge instead of this leftmost thorn.
\item[(ii)] If the leftmost child is an edge connected to black vertex $v$, then $\tau$ is equivalent to the pair $(\tau_1,\tau_2)$ in $\mathcal{W} \times \mathcal{B}$ where $\tau_2$ is the cactus tree descending from $v$ with $v$ as the root and the edge linking $v$ to the root of $\tau$ replaced by the planted edge. $\tau_1$ is the remaining cactus tree descending from the root of $\tau$ with the edge linking it to $v$ as the planted edge.
\item[(iii)] If the leftmost child is a triangle containing black vertex $v_1$ and grey vertex $v_2$ then $\tau$ is equivalent to the tuple  $(\tau_1,\tau_2, \tau_3, t_w)$ in $\mathcal{W} \times \mathcal{B} \times \mathcal{G} \times$${TT}_w$ (${TT}_w$ is the singleton set composed of the triangle child of a white vertex) where $\tau_2$ and $\tau_3$ are the descending trees from $v_1$ and $v_2$ with the edge linking $\tau$'s root and $v_1$ and the edge linking $v_1$ and $v_2$ replaced by a planted edge. $\tau_1$ is the remaining descending cactus tree from its root with the leftmost triangle replaced by the planted edge.
\end{compactenum}
One can check easily that the numbers of triangles, white, black, and grey vertices and their degree distribution are stable by the bijective transformation described above. The complicated case above is case (iii) where the edges linking $v_1$ and $v_2$, and the edge linking the white root of  $\tau$ and $v_2$ are replaced by nothing in (respectively) $\tau_2$ and $\tau_1$. However in Definition \ref{degdef} of the degree of a vertex in $\tau$, these edges were already not taken into account for the degree of respectively $v_1$ and the root vertex. As a consequence : 
\begin{equation}
F = W + W\cdot B + W\cdot B\cdot G \cdot T_w
\end{equation}

\begin{figure}[h]
\begin{center}
\includegraphics[width=125mm]{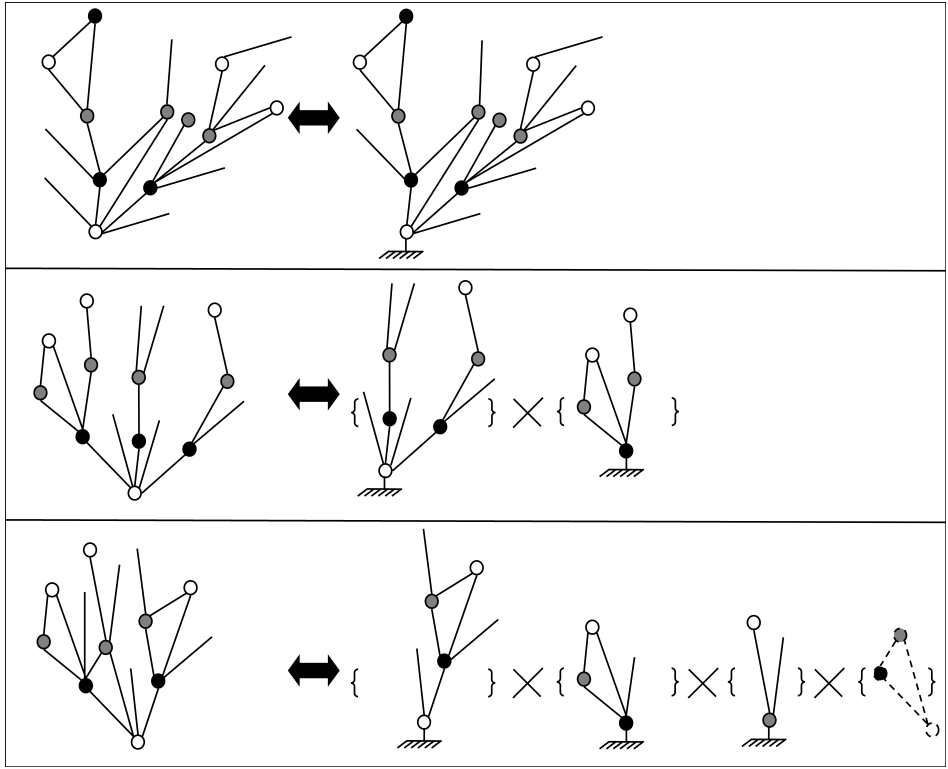}
\caption{Illustration of the decomposition into planted trees}
\label{fig : FW}
\end{center} 
\end{figure}

This decomposition is illustrated in Figure \ref{fig : FW}.

\begin{figure}[h]
\begin{center}
\includegraphics[width=125mm]{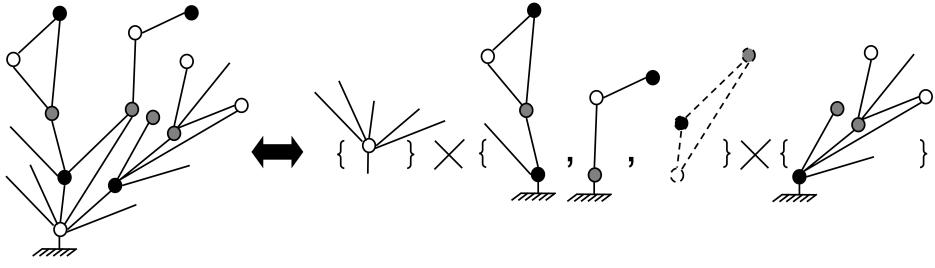}
\caption{Decomposition of a white rooted planted cactus tree}
\label{fig : decomp}
\end{center} 
\end{figure}

To determine F, we show that $W, B, G, T_g, T_w,$ and $T_b$ satisfy a system of functional equations. Namely, as shown in Figure \ref{fig : decomp} any planted cactus tree in $\mathcal{W}$, $\tau$ can be decomposed into:
\begin{compactitem}
\item its white root,
\item the cactus trees rooted in a black vertex descending from the root,
\item a triple composed of a black rooted cactus tree, a grey rooted cactus tree, a triangle for each triangle descending  from the root,
\item the positions of the triangles in the list of children.
\end{compactitem}

\noindent Let $i$ denote the degree of the root vertex (of degree $i+1$), $j$ the number of black children not belonging to a triangle and $k$ the number of descending triangles. The vectors $\mathbf{j}$ and $\mathbf{k}$ give the positions of the $j$ black vertices and $k$ triangles within the $i$ children.  Using the decomposition above, we have:
\begin{equation}
W = x_1\sum_{i \geq 0}t_{i+1}\sum_{0 \leq j+k \leq i}\sum_{{\bf j},{\bf k}}B^{j}(B\cdot G\cdot T_w)^{k}
\end{equation}
\noindent Then:
\begin{eqnarray}
\nonumber W &=& x_1\sum_{i \geq 0}t_{i+1}\sum_{0 \leq j+k \leq i}\binom{i}{j, k}B^j(B\cdot G\cdot T_w)^k\\
\nonumber W &=& x_1\sum_{i \geq 0}t_{i+1}\left (1 + B + B\cdot G\cdot T_w \right)^i%\\
%\label{eq W}  &=& x_1\Phi_1(W,B,G,T_g,T_w,T_b)
\end{eqnarray}
\noindent Similarly,
\begin{eqnarray}
\nonumber B &=& x_2 \sum_{i \geq 0}u_{i+1}\left (1 + G + G\cdot W\cdot T_b \right)^i \\% \label{eq B} &=& x_2\Phi_2(W,B,G,T_g,T_w,T_b)\\
\nonumber G &=& x_3\sum_{i \geq 0}v_{i+1}\left (1 + W + W\cdot B\cdot T_g \right)^i%\\ \label{eq G} &=& x_3\Phi_3(W,B,G,T_g,T_w,T_b)
\end{eqnarray}
Finally :
\begin{equation}
(W,B,G,T_g,T_w,T_b) = {\bf x}{\bf \Phi}(W,B,G,T_g,T_w,T_b)
\end{equation}
\noindent where ${\bf x} = (x_1,x_2,x_3,x_4,x_5,x_6)$ and ${\bf \Phi} = (\Phi_i)_{1\leq i \leq 6}$ with:
\begin{align} 
&\Phi_1(W,B,G,T_g,T_w,T_b) = \sum_{i \geq 0}t_{i+1}\left (1 + B + B\cdot G\cdot T_w \right)^i\\
&\Phi_2(W,B,G,T_g,T_w,T_b) = \sum_{i \geq 0}u_{i+1}\left (1 + G + G\cdot W\cdot T_b \right)^i\\
&\Phi_3(W,B,G,T_g,T_w,T_b) = \sum_{i \geq 0}v_{i+1}\left (1 + W + W\cdot B\cdot T_g \right)^i\\
&\Phi_i = 1\mbox{ for }4\leq i \leq 6
\end{align}
\noindent Using the multivariate Lagrange inversion formula for monomials (see \cite[1.2.9]{GJCE}), we find :
\begin{eqnarray}
\nonumber {k_1 k_2 k_3 k_4 k_5 k_6}\,[{\bf x}^{\bf k}]\,\,W^{r_1}B^{r_2}G^{r_3}T_w^{r_5} =\hspace{0mm}&& \\
\label{eq lagrange}\sum_{\{\mu_{ij}\}}\mid\mid \delta_{ij}k_j-\mu_{ij}\mid\mid\prod_{1\leq i \leq 6}&&\hspace{-8mm}[W^{\mu_{i1}}B^{\mu_{i2}}G^{\mu_{i3}}T_g^{\mu_{i4}}T_w^{\mu_{i5}}T_b^{\mu_{i6}}]\Phi_i^{k_i}
\end{eqnarray}
\noindent where $\mid\mid \cdot \mid\mid$ denotes the determinant, ${\bf k} = (k_1,k_2,k_3,k_4,k_5,k_6)$, $\delta_{ij}$ is the Kronecker delta function and the sum is over all $6\times 6$ integer matrices $\{\mu_{ij}\}$ such that:

\begin{minipage}{0.4\linewidth}\centering
\begin{compactitem}
\item $\mu_{11} = \mu_{14} = \mu_{16} = \mu_{22} = \mu_{24}=\mu_{25} = \mu_{33} =\mu_{35}=\mu_{36} = 0$
\item $\mu_{ij} = 0$ for $i\geq 4$
\item $\mu_{21} + \mu_{31} = k_1 - r_1$
\item $\mu_{12} + \mu_{32} = k_2-r_2$
\item $\mu_{13} + \mu_{23} = k_3-r_3$
\item$\mu_{34} = k_4$
\item $\mu_{15} = k_5-r_5$
\item $\mu_{26} = k_6$
\end{compactitem}
\end{minipage}
\begin{minipage}{0.5\linewidth}
\centering
{i.e.} $\mu=
\left[
\begin{array}{cccccc}
0 & * & * & 0& * & 0\\
* & 0 & *&0&0&*\\
* &* & 0 & * & 0 & 0\\
0 & 0 & 0 & 0& 0 & 0\\
0 & 0 & 0 & 0& 0 & 0\\
0 & 0 & 0 & 0& 0 & 0
\end{array}
\right]
$
\end{minipage}

Looking for zero contribution terms in expression \eqref{eq lagrange}, we notice that $G$ and $T_w$ have necessarily the same degree in the formal power series expansion of $\Phi_1$. Hence, a non zero contribution of 
$$[W^{\mu_{11}}B^{\mu_{12}}G^{\mu_{13}}T_1^{\mu_{14}}T_2^{\mu_{15}}T_3^{\mu_{16}}]\,\,\Phi_1^{k_1}$$ 
implies $\mu_{13} = \mu_{15} = k_5-r_5$. Similar remarks give non zero contributions only for $\mu_{21} = k_6$ and $\mu_{32} = k_4$. As a result, only that one matrix $\mu$ yields a non zero contribution.\\
\noindent For this particular $\mu$, 
\begin{eqnarray}
\nonumber \frac{1}{k_4 k_5 k_6}\mid\mid \delta_{ij}k_j-\mu_{ij}\mid\mid &=& r_1 \left ( k_2 k_3 - (k_3-k_5)k_4 \right )\\
\nonumber &+&r_2 \left( k_6 k_3 + (k_6+r_1-k_1)(r_3+k_4-r_5-k_3) \right) \\
&+& r_3   \left (k_4 k_6-k_2(k_6+r_1-k_1)  \right)\\
&=& \frac{1}{k_4 k_5 k_6} \Delta(\mathbf{k},\mathbf{r}).
\end{eqnarray}
Let $co(k)$ denotes the set of sequences of non negative integers of total sum $k$. The next step is to notice that:
\begin{eqnarray}
\nonumber \Phi_1^{k_1} &=&\sum_{\mathbf{s}  \in co(k_1)}\binom{k_1}{\mathbf{s}}\prod_{i \geq 0}\left [t_{i+1}\left (1 + B + B\cdot G\cdot T_w \right)^i \right]^{s_i} \\
\Phi_1^{k_1} &=& \sum_{\mathbf{s}  \in co(k_1)}\binom{k_1}{\mathbf{s}}\prod_{i \geq 0}t_{i+1}^{s_i} \sum_{a_1,a_2} \binom{\sum_{i} i s_i}{a_1-a_2,a_2}B^a_1(GT_w)^{a_2}.
\end{eqnarray}
As a result, the coefficient in $W^{\mu_{11}}B^{\mu_{12}}G^{\mu_{13}}T_g^{\mu_{14}}T_w^{\mu_{15}}T_b^{\mu_{16}}$ is equal to
\begin{equation}
\sum_{\mathbf{s}  \in co(k_1)}\binom{k_1}{\mathbf{s}}\prod_{i \geq 0}t_{i+1}^{s_i} \binom{\sum_{i} i s_i}{k_5-r_5,k_2-k_4-k_5-r_2+r_5}.
\end{equation}
Similarly, we have:
\begin{eqnarray}
\nonumber [W^{\mu_{21}}B^{\mu_{22}}G^{\mu_{23}}T_g^{\mu_{24}}T_w^{\mu_{25}}T_b^{\mu_{26}}]\,\,\Phi_2^{k_2} =&&\\
\sum_{\mathbf{s}  \in co(k_2)}\binom{k_2}{\mathbf{s}}\prod_{i \geq 0}u_{i+1}^{s_i}&&\hspace{-7mm} \binom{\sum_{i} i s_i}{k_6,k_3-k_5-k_6-r_3+r_5}\\
\nonumber \left [W^{\mu_{31}}B^{\mu_{32}}G^{\mu_{33}}T_g^{\mu_{34}}T_w^{\mu_{35}}T_b^{\mu_{36}}\right ]\,\,\Phi_3^{k_3} =&&\\
\sum_{\mathbf{s}  \in co(k_3)}\binom{k_3}{\mathbf{s}}\prod_{i \geq 0}v_{i+1}^{s_i}&&\hspace{-7mm} \binom{\sum_{i} i s_i}{k_4,k_1-k_4-k_6- r_1}.
\end{eqnarray}
Putting everything together gives:
\begin{eqnarray}
\nonumber  [x_1^{\ell(\lambda)}x_2^{\ell(\mu)}x_3^{\ell(\nu)}x_4^{g}x_5^{w}x_6^{b}\mathbf{t}^{\mathbf{n}(\lambda)}\mathbf{u}^{\mathbf{n}(\mu)}\mathbf{v}^{\mathbf{n}(\nu)}]\,\,&&\hspace{-7mm}W^{r_1}B^{r_2}G^{r_3}T_w^{r_5}= \\
\nonumber \\
\nonumber \frac{ \Delta(\ell(\lambda),\ell(\mu),\ell(\nu),g,w,b,\mathbf{r})}{\ell(\lambda)\ell(\mu) \ell(\nu)}&&\hspace{-7mm} \binom{\ell(\lambda)}{\mathbf{n}(\lambda)}\binom{\sum_{i} i n_{i+1}(\lambda)}{w-r_5,\ell(\mu)-g-w-r_2+r_5} \times \\ 
\nonumber \times &&\hspace  {-7mm} \binom{\ell(\mu)}{\mathbf{n}(\mu)}\binom{\sum_{i} in_{i+1}(\mu)}{b,\ell(\nu)-w-b-r_3+r_5}\\
\times&&\hspace{-7mm}\binom{\ell(\nu)}{\mathbf{n}(\nu)} \binom{\sum_{i} i n_{i+1}(\nu)}{g,\ell(\lambda)-g-b- r_1}.
\end{eqnarray}
Noticing that for $\epsilon \vdash n$
\begin{equation}
\sum_{i \geq 0}in_{i+1}(\epsilon) = \sum_{i \geq 0}(i+1)n_{i+1}(\epsilon) - \sum_{i \geq 0}n_{i+1}(\epsilon) = n-\ell(\epsilon).
\end{equation} 
\noindent And summing for $\mathbf{r} \in \{(1,0,0,0), (1,1,0,0), (1,1,1,1)\}$ gives the desired result.

\section{Proof of Corollary \ref{cor1} and restriction of bijection $\Theta$ when $\nu=[1^n]$} \label{sect:pfcor1}

\noindent We look more closely at the case when one of the partitions, say $\nu$, is $[1^n]$. We need the following definitions:

\begin{defn}[Partitioned bicolored map]
 Given partitions $\la, \mu \vdash n$, let $\mathcal{C}(\la, \mu)$ be the set of triples $(\pi_1,\pi_2,\al)$ such that $\al \in \mathfrak{S}_n$, $\pi_1,\pi_2$ are set partitions of $[n]$ with $\ty(\pi_1) =\la$ and $\ty(\pi_2) = \mu$, and each block of $\pi_1$ and $\pi_2$ is a union of cycles of $\al$ and $\be = \al^{-1}\ga_n$ respectively. The elements of $\mathcal{C}(\la,\mu)$ are called {\em unicellular partitioned bicolored maps} of type $\la$ and $\mu$. Let $C(\la,\mu) = |\mathcal{C}(\la,\mu)|$.
\end{defn}

\begin{defn}[Ordered rooted bicolored thorn trees]
For $\la, \mu \vdash n$ such that $\ell(\la)+\ell(\mu) \leq n+1$, we define $\widetilde{\mathcal{BT}}(\la, \mu)$ as the set of ordered rooted bicolored trees with $\ell(\la)$  white vertices, $\ell(\mu)$  black vertices, $n+1-\ell(\la)-\ell(\mu)$ thorns connected to the black vertices  and $n+1-\ell(\la)-\ell(\mu)$  thorns connected to the white vertices. The white (respectively black) vertices' degree distribution (accounting the thorns) is specified by $\la$ (respectively $\mu$). The root is a white vertex. 
\end{defn} 

Again, adapting the Lagrange inversion developed in \cite{GJ92}, we get:
$$
|\widetilde{BT}(\lambda, \mu)| = \frac{n}{Aut(\lambda)Aut(\mu)}\frac{(n-\ell(\lambda))!(n-\ell(\mu))!}{{(n+1-\ell(\lambda)-\ell(\mu))!\,}^2}.
$$
We now prove Corollary \ref{cor1}: 
\begin{proof}
We have $\mathcal{C}(\la,\mu,[1^n])=\mathcal{C}(\la,\mu)$, the number of unicellular partitioned bicolored maps of type $\la$ and $\mu$. Indeed, as the cycles of $\al_3$ refine the blocks of $\pi_3$, if $\nu=[1^n]$ then $\pi_3=\{\{1\}, \{2\},\ldots, \{n\}\}$ and $\al_3=\iota$, the identity permutation. Then extracting the coefficient of $m_{1^n}({\bf z})$ to both sides of \eqref{ngs} we obtain
\footnotesize
\begin{eqnarray*}
\sum_{\la,\mu \vdash n}Aut(\la)Aut(\mu)Aut(1^n) C(\la,\mu)m_{\la}({\bf x})m_{\mu}({\bf y}) &=& [m_{1^n}({\bf z})] \sum_{\la,\mu,\nu \vdash n} {k_{\la,\mu,\nu}^n} p_{\la}({\bf x})p_{\mu}({\bf y})p_{\nu}({\bf z})\\
&=& \sum_{\nu \vdash n, \nu \preceq 1^n} Aut(1^n)\overline{R}_{\nu,1^n}\sum_{\la,\mu \vdash n} k_{\la,\mu,\nu}^n p_{\la}({\bf x})p_{\mu}({\bf y})
\end{eqnarray*}
\normalsize
Since $\overline{R}_{\nu,1^n}=1$ if $\nu=1^n$ and zero otherwise, we  obtain 
\begin{eqnarray*}
\sum_{\la,\mu \vdash n} Aut(\la)Aut(\mu)C(\la,\mu)m_{\la}({\bf x})m_{\mu}({\bf y}) &=& \sum_{\la,\mu \vdash n} k_{\la,\mu,1^n}^n p_{\la}({\bf x})p_{\mu}({\bf y}), 
\end{eqnarray*}
where $k_{\la,\mu,1^n}^n = k_{\la,\mu}^n$. 
%By Theorem~\ref{thm1}, the LHS of Equation \eqref{eq2} is $1/Aut(1^n)=1/n!$ times the coefficient of $m_{1^n}({\bf z})$ of the RHS in Equation \eqref{eq1}. Since the parameter $g\leq n-\ell(\nu)=0$, then after some straightforward cancellations we get the desired RHS.
\end{proof}

\noindent Next, we say what the bijection $\Theta^n_{\lambda,\mu,\nu}$ of Theorem \ref{thm:bij} does in this case ($\nu=[1^n]$). This matches the bijection in \cite{MV} which in turn matches the bijection in \cite{GJ92} when ${\sf g}(\la,\mu)=0$ and is a refinement of a bijection in \cite{SV}. 

\begin{cor}
There is a bijection between partitioned bicolored maps $\mathcal{C}(\la,\mu,n)$ and pairs $(\tilde{t},\sigma)$ where $\tilde{t}\in \widetilde{BT}(\la,\mu)$ and $\sigma \in \mathfrak{S}_{n+1-\ell(\la)-\ell(\mu)}$. 
\end{cor}

\begin{proof}
From above we have that $\mathcal{C}(\lambda,\mu,[1^n])=\mathcal{C}(\lambda,\mu)$. Let $$(\widetilde{\tau},\sigma_1,\sigma_2,\chi) :=\Theta^n_{\lambda,\mu,[1^n]}(\pi_1,\pi_2,[1^n],\al_1,\al_1^{-1}\gamma)$$ for $(\pi_1,\pi_2,\al_1)\in\mathcal{C}(\lambda,\mu)$. We know that $\ell(\mu)=n$ forces $\alpha_3=\iota$, the identity permutation. In this case ${m'}_2^{(j)}$ is the maximal element of $\al_2^{-1}(\pi_2^{(j)})$. But $\al_2$ preserves the blocks of $\pi_2$, thus  ${m'}_2^{(j)}$ is just the maximal element of $\pi_2^{(j)}$, call this $m_2^{(j)}$. First, we show that in this case $\widetilde{\tau}$ can be reduced to a tree $\tilde{t} \in \widetilde{\mathcal{BT}}(\la, \mu)$. Then, we show that $\sigma_1,\sigma_2$ are trivial permutations and that $\chi$ can be regarded as a permutation in $\mathfrak{S}_{n+1-\ell(\la)-\ell(\mu)}$.

From the incidence rules in Table \ref{incidrules}, we see that each black vertex $j$ has $|\pi_2^j|$ children (one for each element of the block). And $\ell(\la)-1$ of the grey vertices have one child (one for each $m_1^{(i)}$, $1\leq i \leq \ell(\la)-1$), the other grey vertices have none. Recall $w,b,g$ count the number of triangles in $\widetilde{\tau}$ children of white, black, and grey vertices respectively. From the rules in Table \eqref{rulestriang} for adding triangles children of the different vertices, we see that $w=\ell(\mu)$. And if a grey vertex has a white child, then these two vertices are part of a triangle child of a black vertex, so $b=\ell(\la)-1$. For triangles children of grey vertices,  if $\al_2(m_2^{(j)})=m_1^{(i)}$ for some $i$ and $j$ ($1\leq i \leq \ell(\la)-1$ and $1\leq j \leq \ell(\mu)$), then $m_1^{(i)}\in\pi_2^{(j)}$ and $m_1^{(i)}\leq m_2^{(j)}$ ($\al_2$ preserves blocks of $\pi_2$). But $\al_1(m_1^{(i)})\leq m_1^{(i)}$   ($\al_1$ preserves blocks of $\pi_1$), so $\gamma(m_2^{(j)})\leq m_2^{(j)}$. This only happens if $m_2^{(j)}=n$ which means $1=\gamma(n)\in \pi^{(i)}$ and $i=\ell(\la)$, a contradiction. Thus $g=0$; there are no triangles children of a grey vertex.

In terms of the thorns, the cactus $\widetilde{\tau}$ has $n+1-\ell(\la)-\ell(\mu)$ thorns connected to white vertices and since $n-\ell(\mu)-\ell(\nu)+w=0$ and $n+1-\ell(\la)-\ell(\nu)+b=0$, $\widetilde{\tau}$ has no thorns connected to black and grey vertices.

From above we see that each grey vertex is either: (i) within a triangle child of a black vertex, (ii) a vertex of a triangle child of a white vertex, and (iii) a leaf (note that there are $n-(\ell(\la)-1)-\ell(\mu)$ of these). Then depending on the case we do the following reductions: (i) and (ii) triangle to the edge linking the white and the black vertex, (iii) leaf to thorn connected to a black vertex. We summarize this reduction graphically in Table \ref{reducrules}:
\begin{table}
$$
\begin{array}{|l|l|l|} \hline
%\multicolumn{3}{|l|}{\text{Local rules for reducing cactus tree $\widetilde{\tau}$ when $\nu=[1^n]$}} \\ \hline
\begin{array}{ccc}
\begin{tikzpicture}[thick,scale=0.8]
\node (r) at ( 0,0) [bv] {}; 
\node (gc) at (0,1) [gv] {};
\node (gc2) at (0.5,2) [wv] {};
\draw [-] (r) -- (gc); 
\draw [-] (r) -- (gc2);
\draw [-] (gc) -- (gc2); 
\end{tikzpicture}
&
\begin{tikzpicture}
\draw[->]	(0,0.8)	-- (0.5,0.8);
\node (1) at (0,0) {};
\end{tikzpicture}
&
\begin{tikzpicture}[thick,scale=0.8]
\node (r) at ( 0,0) [bv] {}; 
\node (gc2) at (0.5,2) [wv] {}; 
\draw [-] (r) -- (gc2);
\end{tikzpicture}
\end{array}
&
\begin{array}{ccc}
\begin{tikzpicture}[thick,scale=0.8]
\node (r) at ( 0,0) [bv] {}; 
\node (gc) at (0,1) [gv] {};
\node (gc2) at (0.5,2)  {};
\draw [-] (r) -- (gc); 
\end{tikzpicture}
&
\begin{tikzpicture}
\draw[->]	(0,0.8)	-- (0.5,0.8);
\node (1) at (0,0) {};
\end{tikzpicture}
&
\begin{tikzpicture}[thick,scale=0.8]
\node (r) at ( 0,0) [bv] {}; 
\node (gc) at (0,1) {};
\node (gc2) at (0.5,2)  {};
\draw [-] (r) -- (gc); 
\end{tikzpicture}
\end{array}
&
\begin{array}{ccc}
\begin{tikzpicture}[thick,scale=0.8]
\node (r) at ( 0,0) [wv] {}; 
\node (gc) at (0,1) [bv] {};
\node (gc2) at (0.5,2) [gv] {};
\draw [-] (r) -- (gc); 
\draw [-] (r) -- (gc2);
\draw [-] (gc) -- (gc2); 
\end{tikzpicture}
&
\begin{tikzpicture}
\draw[->]	(0,0.8)	-- (0.5,0.8);
\node (1) at (0,0) {};
\end{tikzpicture}
&
\begin{tikzpicture}[thick,scale=0.8]
\node (r) at ( 0,0) [wv] {}; 
\node (gc) at (0,1) [bv] {};
\draw [-] (r) -- (gc);  
\end{tikzpicture}
\end{array} \\ \hline
\end{array}
$$
\caption{Local rules for reducing cactus tree $\widetilde{\tau}$ when $\nu=[1^n]$.}
\label{reducrules}
\end{table}

The outcome is an ordered bicolored tree $\tilde{t}$ with $\ell(\la)$ white vertices and $\ell(\mu)$ black vertices. This tree $\tilde{t}$ has $n+1-\ell(\la)-\ell(\mu)$ thorns connected to white vertices and $n+1-\ell(\la)-\ell(\mu)$ thorns connected to black vertices. Moreover, this reduction $\widetilde{\tau} \to \tilde{t}$ is reversible.

In addition, since $\widetilde{\tau}$ had no thorns connected to black and grey vertices ($n+1-\ell(\la)-\ell(\nu)+b=0$ and $n-\ell(\mu)-\ell(\nu)+w=0$), then $\sigma_1$ and $\sigma_2$ are trivial permutations. Since $\ell(\nu)-w-b=n+1-\ell(\la)-\ell(\mu)=n+1-\ell(\la)-\ell(\mu)+g$, then we see that $\chi$ is just a permutation $\sigma$ in $\mathfrak{S}_{n+1-\ell(\la)-\ell(\mu)}$.

In summary, we have a bijection from $\mathcal(\la,\mu,n)$ to the desired pair $(\tilde{t},\sigma)$.
\end{proof}

\begin{exm}
Let $\alpha_1=(189\,10)(25)(3467)$, $\alpha_2=(15427)(3)(6)(8)(9)(10)$, $\alpha_3=\iota$ ($\al_1\al_2=\al_1\al_2\al_3=\gamma_{10}$), $\pi_1=\left\{\{3,4,6,7\},\{1,2,5,8,9,10\}\right\}$, $\pi_2 = \left\{\{1,2,4,5,7,10\},\{3,9\},\{6,8\} \right\}$, $\pi_3=\left\{\{1\},\{2\},\ldots,\{10\}\right\}$. Then $\Theta_{\lambda,\mu,[1^n]}^n(\pi_1,\pi_2,\pi_3,\alpha_1,\alpha_2)=(\widetilde{\tau},\emptyset,\emptyset,251364)$ where $\widetilde{\tau}$ and its reduction $\tilde{t}$ are depicted below: 

\begin{figure}[htbp]
$$
\begin{array}{ccc}
\begin{tikzpicture}
                \node [style=wv] (0) at (1.5, 3) {};
		\node [style=gv] (1) at (-2, 2) {};
		\node [style=gv] (2) at (-1, 2) {};
		\node [style=gv] (3) at (-0.5, 2) {};
		\node [style=gv] (4) at (0.5, 2) {};
		\node [style=gv] (5) at (1, 2) {};
		\node [style=gv] (6) at (1.5, 2) {};
		\node [style=gv] (7) at (2, 2) {};
		\node [style=gv] (8) at (2.6, 2) {};
		\node [style=gv] (9) at (3.2, 2) {};
		\node [style=gv] (10) at (3.9, 2) {};
		\node [style=bv] (11) at (-1, 1) {};
		\node [style=bv] (12) at (0, 1) {};
		\node [style=bv] (13) at (1, 1) {};
		\node [style=wv] (14) at (0, 0) {};
                \node (t1) at (1.5,0.5) {};
                \node (t2) at (2,0.5) {};
                \node (t3) at (2.5,0.5) {};
                \node (t4) at (1.1,4) {};
                \node (t5) at (1.5,4) {};
                \node (t6) at (1.9,4) {};
                \draw [-] (14) -- (t1);
                \draw [-] (14) -- (t2);
                \draw [-] (14) -- (t3);
                \draw [-] (0) -- (t4);
                \draw [-] (0) -- (t5);
                \draw [-] (0) -- (t6);
                \draw [-] (14) -- (13);
                \draw [-] (14) -- (12);
                \draw [-] (14) -- (11);
                \draw [-] (13) -- (10);
                \draw [-] (13) -- (9);
                \draw [-] (13) -- (8);
                \draw [-] (13) -- (7);
                \draw [-] (13) -- (6);
                \draw [-] (13) -- (5);
                \draw [-] (13) -- (0);
                \draw [-] (12) -- (4);
                \draw [-] (12) -- (3);
                \draw [-] (11) -- (2);
                \draw [-] (11) -- (1);
                \draw [-] (6) -- (0);
                \draw [-] (14) -- (2);
                \draw [-] (14) -- (4);
                \draw [-] (14) -- (10);
\end{tikzpicture}
& 
\begin{tikzpicture}
\draw[->]	(0,2)	-- (1,2);
\node (1) at (0,0) {};
\end{tikzpicture}
&
\begin{tikzpicture}
                \node [style=wv] (0) at (1.5, 3) {};
		\node (1) at (-2, 2) {};
		\node (3) at (-0.5, 2) {};
		\node (5) at (1, 2) {};
		\node (7) at (2, 2) {};
		\node (8) at (2.6, 2) {};
		\node (9) at (3.2, 2) {};
		\node [style=bv] (11) at (-1, 1) {};
		\node [style=bv] (12) at (0, 1) {};
		\node [style=bv] (13) at (1, 1) {};
		\node [style=wv] (14) at (0, 0) {};
                \node (t1) at (1.5,0.5) {};
                \node (t2) at (2,0.5) {};
                \node (t3) at (2.5,0.5) {};
                \node (t4) at (1.1,4) {};
                \node (t5) at (1.5,4) {};
                \node (t6) at (1.9,4) {};
                \draw [-] (14) -- (t1);
                \draw [-] (14) -- (t2);
                \draw [-] (14) -- (t3);
                \draw [-] (0) -- (t4);
                \draw [-] (0) -- (t5);
                \draw [-] (0) -- (t6);
                \draw [-] (14) -- (13);
                \draw [-] (14) -- (12);
                \draw [-] (14) -- (11);
                \draw [-] (13) -- (9);
                \draw [-] (13) -- (8);
                \draw [-] (13) -- (7);
                \draw [-] (13) -- (5);
                \draw [-] (13) -- (0);
                \draw [-] (12) -- (3);
                \draw [-] (11) -- (1);
\end{tikzpicture}
\end{array}
$$
\caption{Example of reduction of a cactus tree to a bicolored thorn trees when $\nu=[1^{10}]$.} 
\end{figure}
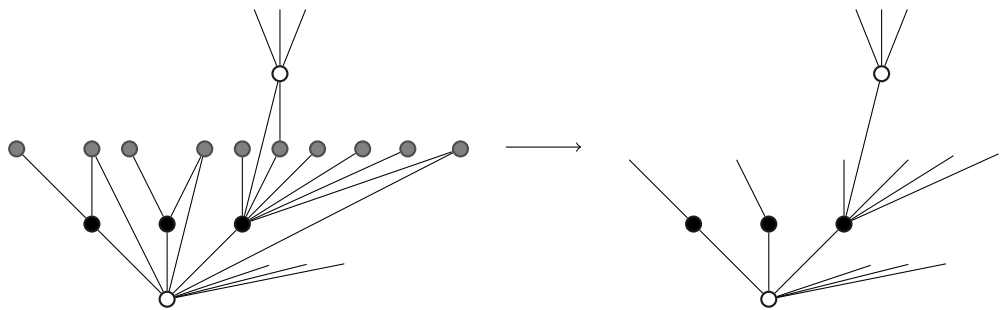
\end{exm}

\section*{Acknowledgements}
 
We thank the anonymous referee for several helpful comments and suggestions. A.H.M. thanks the support and hospitality of the Combinatorial Models Group of \'Ecole Polytechnique where part of this work was done.

\small

%Department of Mathematics, Massachusetts Institute of Technology \\
%Cambridge, MA USA 02139\\
%{\tt ahmorales@math.mit.edu}

\end{document}